# STOCHASTIC CALCULUS FOR SYMMETRIC MARKOV PROCESSES

By Z.-Q. Chen,[1] P. J. Fitzsimmons,[2] K. Kuwae[3] and T.-S. Zhang[4]

*University of Washington, University of California at San Diego, Kumamoto University and University of Manchester*

*Dedicated to S. Nakao on the occasion of his 60th birthday*

Using time-reversal, we introduce a stochastic integral for zero-energy additive functionals of symmetric Markov processes, extending earlier work of S. Nakao. Various properties of such stochastic integrals are discussed and an Itô formula for Dirichlet processes is obtained.

**1. Introduction and framework.** It is well known that stochastic integrals and Itô's formula for semimartingales play a central role in modern probability theory. However, there are many important classes of Markov processes that are not semimartingales. For example, symmetric diffusions on $\mathbb{R}^d$ whose infinitesimal generators are elliptic operators in divergence form $\mathcal{L} = \sum_{i,j=1}^{d} \frac{\partial}{\partial x_i}(a_{ij}(x)\frac{\partial}{\partial x_j})$ with merely measurable coefficients need not be semimartingales. Even when $X$ is a Brownian motion in $\mathbb{R}^d$ and $u \in W^{1,2}(\mathbb{R}^d) := \{u \in L^2(\mathbb{R}^d; dx) \mid |\nabla u| \in L^2(\mathbb{R}^d; dx)\}$, the process $u(X_t)$ is not generally a semimartingale. To study such processes, Fukushima ob-

Received August 2006; revised April 2007.
[1]Supported in part by NSF Grant DMS-06-00206.
[2]Supported by a foundation based on the academic cooperation between Yokohama City University and UCSD.
[3]Supported by a foundation based on the academic cooperation between Yokohama City University and UCSD, and partially supported by a Grant-in-Aid for Scientific Research (C) No. 16540201 from the Japan Society for the Promotion of Science.
[4]Supported in part by the British EPSRC.

tained the following substitute for Itô's formula (see [7]): for $u \in W^{1,2}(\mathbb{R}^d)$,

(1.1) $$u(X_t) = u(X_0) + M_t^u + N_t^u \qquad \text{for } t \geq 0,$$

$\mathbf{P}_x$-a.s. for quasi-every $x \in \mathbb{R}^d$, where $M^u$ is a square-integrable martingale and $N^u$ is a continuous additive functional of zero energy. The decomposition (1.1) is called *Fukushima's decomposition* and holds for a general symmetric Markov process $X$ and for $u \in \mathcal{F}$, where $(\mathcal{E}, \mathcal{F})$ is the Dirichlet space for $X$. In this paper, a stochastic process $\xi = \{\xi_t, t \geq 0\}$ under some $\sigma$-finite measure $\mathbf{P}$ is called a *Dirichlet process* if $\xi$ has locally finite quadratic variation under $\mathbf{P}$. The composite process $u(X)$ is a Dirichlet process under $\mathbf{P}_m$, where $m$ is the Lebesgue measure on $\mathbb{R}^d$, as it has finite quadratic variation on compact time intervals. Nakao introduced a stochastic integral $\int_0^t f(X_s) \, dN_s^u$ in [14] by using a Riesz representation theorem in a suitably constructed Hilbert space. Nakao's stochastic integral played an important role in the study of lower order perturbation of diffusion processes by Lunt, Lyons and Zhang [12] and by Fitzsimmons and Kuwae [5]. However, Nakao's definition of the stochastic integral $\int_0^t f(X_s) \, dN^u$, requiring $u$ to be in the domain of the Dirichlet form of $X$ and $f$ to be square-integrable with respect to the energy measure of $u$, is too restrictive to be useful in the study of lower-order perturbation for symmetric Markov processes with discontinuous sample paths, such as stable processes. Such a study requires stochastic integrals for more general integrators as well as integrands. The purpose of this paper is to present a new way of defining the stochastic integral for Dirichlet processes associated with a symmetric Markov process. Our new approach uses only the time-reversal operator for the process $X$ and is therefore more direct and provides additional insight into stochastic integration for Dirichlet processes. This approach enables us to define $\Lambda(M)$ [see (1.5)] for any locally square-integrable martingale additive functional (MAF) $M$, subject to some mild conditions. Thus, it not only recovers Nakao's results in [14], but also extends them significantly. The new stochastic integral allows us to study various transforms for symmetric Markov processes, a project that is carried out in a subsequent paper [2].

A more detailed description of the current paper appears below.

Let $X = \{\Omega, \mathcal{F}_\infty, \mathcal{F}_t, X_t, \theta_t, \zeta, \mathbf{P}_x, x \in E\}$ be an $m$-symmetric right Markov process with a Lusin state space $E$, where $m$ is a $\sigma$-finite measure with full support on $E$. Its associated Dirichlet space $(\mathcal{E}, \mathcal{F})$ on $L^2(E; m)$ is known to be quasi-regular (see [13]). By [1], $(\mathcal{E}, \mathcal{F})$ is quasi-homeomorphic to a regular Dirichlet space on a locally compact separable metric space. Using this quasi-homeomorphism, there is no loss of generality in assuming that $X$ is an $m$-symmetric Hunt process on a locally compact metric space $E$ such that its associated Dirichlet space $(\mathcal{E}, \mathcal{F})$ is regular on $L^2(E; m)$ and that $m$ is a positive Radon measure with full topological support on $E$. We assume this throughout the sequel.

Without loss of generality, we can take $\Omega$ to be the canonical path space $D([0,\infty[ \to E_\Delta)$ of right-continuous, left-limited (*rcll*, for short) functions from $[0,\infty[$ to $E_\Delta$, for which $\Delta$ is a trap [i.e., if $\omega(t) = \Delta$, then $\omega(s) = \Delta$ for all $s > t$]. For any $\omega \in \Omega$, we set $X_t(\omega) := \omega(t)$. Let $\zeta(\omega) := \inf\{t \geq 0 \mid X_t(\omega) = \Delta\}$ be the lifetime of $X$. As usual, $\mathcal{F}_\infty$ and $\mathcal{F}_t$ are the minimal augmented $\sigma$-algebras obtained from $\mathcal{F}^0_\infty := \sigma\{X_s \mid 0 \leq s < \infty\}$ and $\mathcal{F}^0_t := \sigma\{X_s \mid 0 \leq s \leq t\}$, respectively, under $\mathbf{P}_x$; see the next section for more details. We sometimes use a filtration denoted by $(\mathcal{M}_t)$ on $(\Omega, \mathcal{M})$ in order to represent several filtrations, for example, $(\mathcal{F}^0_t)$, $(\mathcal{F}^0_{t+})$ on $(\Omega, \mathcal{F}^0_\infty)$, $(\mathcal{F}_t)$ on $(\Omega, \mathcal{F}_\infty)$ and others introduced later. We use $\theta_t$ to denote the shift operator defined by $\theta_t(\omega)(s) := \omega(t+s)$, $t,s \geq 0$. Let $\omega_\Delta$ be the path starting from $\Delta$. Then, $\omega_\Delta(s) \equiv \Delta$ for all $s \in [0,\infty[$. Note that $\theta_{\zeta(\omega)}(\omega) = \omega_\Delta$ if $\zeta(\omega) < \infty$, $\{\omega_\Delta\} \in \mathcal{F}^0_0 \subset \mathcal{F}^0_t$ for all $t > 0$ and $\mathbf{P}_x(\{\omega_\Delta\}) \leq \mathbf{P}_x(X_0 = \Delta) = 0$ for $x \in E$. For a Borel subset $B$ of $E$, $\tau_B := \inf\{t > 0 \mid X_t \notin B\}$ (the *exit time* of $B$) is an $(\mathcal{F}_t)$-stopping time. If $B$ is closed, then $\tau_B$ is an $(\mathcal{F}^0_{t+})$-stopping time. Also, $\zeta$ is an $(\mathcal{F}^0_t)$-stopping time because $\{\zeta \leq t\} = \{X_t = \Delta\} \in \mathcal{F}^0_t$, $t \geq 0$. The transition semigroup of $X$, $\{P_t, t \geq 0\}$, is defined by

$$P_t f(x) := \mathbf{E}_x[f(X_t)] = \mathbf{E}_x[f(X_t) : t < \zeta], \qquad t \geq 0.$$

Each $P_t$ may be viewed as an operator on $L^2(E;m)$; collectively, these operators form a strongly $L^2$-continuous semigroup of self-adjoint contractions. The Dirichlet form associated with $X$ is the bilinear form

$$\mathcal{E}(u,v) := \lim_{t \downarrow 0} \frac{1}{t}(u - P_t u, v)_m$$

defined on the space

$$\mathcal{F} := \left\{ u \in L^2(E;m) \, \Big| \, \sup_{t>0} t^{-1}(u - P_t u, u)_m < \infty \right\}.$$

Here, we use the notation $(f,g)_m := \int_E f(x)g(x)m(dx)$.

For the reader's convenience, we recall the following definitions from [13] and [7].

DEFINITION 1.1. (i) An increasing sequence $\{F_n\}_{n \geq 1}$ of closed subsets of $E$ is an $\mathcal{E}$-*nest* (or simply *nest*) if and only if $\bigcup_{n \geq 1} \mathcal{F}_{F_n}$ is $\mathcal{E}_1$-dense in $\mathcal{F}$, where $\mathcal{E}_1 = \mathcal{E} + (\cdot, \cdot)_{L^2(E,m)}$ and

$$\mathcal{F}_{F_n} := \{u \in \mathcal{F} : u = 0 \text{ $m$-a.e. on } E \setminus F_n\}.$$

(ii) A subset $N \subset E$ is $\mathcal{E}$-*polar* if and only if there is an $\mathcal{E}$-nest $\{F_n\}_{n \geq 1}$ such that $N \subset \bigcap_{n \geq 1}(E \setminus F_n)$.

(iii) A function $f$ on $E$ is said to be *quasi-continuous* if there is an $\mathcal{E}$-nest $\{F_n\}_{n\geq 1}$ such that $f|_{F_n}$ is continuous on $F_n$ for each $n \geq 1$; we denote this situation briefly by writing $f \in C(\{F_n\})$.

(iv) A statement depending on $x \in A$ is said to hold *quasi-everywhere* (q.e. in abbreviation) on $A$ if there is an $\mathcal{E}$-polar set $N \subset A$ such that the statement is true for every $x \in A \setminus N$.

(v) A nearly Borel subset $N \subset E$ is called *properly exceptional* if $m(N) = 0$ and
$$\mathbf{P}_x(X_t \in E_\Delta \setminus N \text{ for } t \geq 0 \text{ and } X_{t-} \in E_\Delta \setminus N \text{ for } t > 0) = 1$$
for every $x \in E \setminus N$.

It is known (cf. [7]) that a family $\{F_n\}$ of closed sets is an $\mathcal{E}$-nest if and only if
$$\mathbf{P}_x\left(\lim_{n\to\infty} \tau_{F_n} = \zeta\right) = 1 \qquad \text{for q.e. } x \in E.$$

It is also known that a properly exceptional set is $\mathcal{E}$-polar and that every $\mathcal{E}$-polar set is contained in a properly exceptional set. Every element $u$ in $\mathcal{F}$ admits a quasi-continuous $m$-version. We assume throughout this section that functions in $\mathcal{F}$ are always represented by their quasi-continuous $m$-versions. In the sequel, the abbreviations CAF, PCAF and MAF stand for "continuous additive functional," "positive continuous additive functional" and "martingale additive functional," respectively; the definitions of these terms can be found in [7].

Let $\overset{\circ}{\mathcal{M}}$ and $\mathcal{N}_c$ denote, respectively, the space of MAF's of finite energy and the space of continuous additive functionals of zero energy. For $u \in \mathcal{F}$, Fukushima's decomposition holds:

(1.2) $\qquad u(X_t) - u(X_0) = M_t^u + N_t^u \qquad$ for every $t \in [0, \infty[$,

$\mathbf{P}_x$-a.s. for q.e. $x \in E$, where $M^u \in \overset{\circ}{\mathcal{M}}$ and $N^u \in \mathcal{N}_c$.

A positive continuous additive functional (PCAF) of $X$ (call it $A$) determines a measure $\mu = \mu_A$ on the Borel subsets of $E$ via the formula

(1.3) $\qquad \mu(f) = \uparrow \lim_{t\downarrow 0} \frac{1}{t} \mathbf{E}_m\left[\int_0^t f(X_s)\, dA_s\right],$

in which $f: E \to [0, \infty]$ is Borel measurable. Here, $\uparrow \lim_{t\downarrow 0}$ indicates an increasing limit as $t \downarrow 0$. The measure $\mu$ is necessarily *smooth*, in the sense that $\mu$ charges no $\mathcal{E}$-polar set of $X$ and there is an $\mathcal{E}$-nest $\{F_n\}$ of closed subsets of $E$ such that $\mu(F_n) < \infty$ for each $n \in \mathbb{N}$. Conversely, given a smooth measure $\mu$, with $A = A^\mu$. In the sequel, we refer to this bijection between smooth measures and PCAF's as the *Revuz correspondence* and to $\mu$ as the *Revuz measure* of $A^\mu$.

If $M$ is a locally square-integrable martingale additive functional (MAF) of $X$ on the random time interval $[\![0, \zeta[\![$, then the process $\langle M \rangle$ (the dual predictable projection of $[M]$) is a PCAF (Proposition 2.8) and the associated Revuz measure [as in (1.3)] is denoted by $\mu_{\langle M \rangle}$. More generally, if $M^u$ is the martingale part in the Fukushima decomposition of $u \in \mathcal{F}$, then $\langle M^u, M \rangle$ is a CAF locally of bounded variation and we have the associated Revuz measure $\mu_{\langle M^u, M \rangle}$, which is locally the difference of smooth (positive) measures. For $u \in \mathcal{F}$, the Revuz measure $\mu_{\langle M^u \rangle}$ of $\langle M^u \rangle$ will usually be denoted by $\mu_{\langle u \rangle}$.

Let $(N(x, dy), H_t)$ be a Lévy system for $X$; that is, $N(x, dy)$ is a kernel on $(E_\Delta, \mathcal{B}(E_\Delta))$ and $H_t$ is a PCAF with bounded 1-potential such that for any nonnegative Borel function $\phi$ on $E_\Delta \times E_\Delta$ vanishing on the diagonal and any $x \in E_\Delta$,

$$\mathbf{E}_x\left[\sum_{s \leq t} \phi(X_{s-}, X_s)\right] = \mathbf{E}_x\left[\int_0^t \int_{E_\Delta} \phi(X_s, y) N(X_s, dy) \, dH_s\right].$$

To simplify notation, we will write

$$N\phi(x) := \int_{E_\Delta} \phi(x, y) N(x, dy)$$

and

$$(N\phi * H)_t := \int_0^t N\phi(X_s) \, dH_s.$$

Let $\mu_H$ be the Revuz measure of the PCAF $H$. The jumping measure $J$ and the killing measure $\kappa$ of $X$ are then given by

$$J(dx, dy) = \tfrac{1}{2} N(x, dy) \mu_H(dx) \quad \text{and} \quad \kappa(dx) = N(x, \{\Delta\}) \mu_H(dx).$$

These measures feature in the Beurling–Deny decomposition of $\mathcal{E}$: for $f, g \in \mathcal{F}$,

$$\mathcal{E}(f, g) = \mathcal{E}^{(c)}(f, g) + \int_{E \times E} (f(x) - f(y))(g(x) - g(y)) J(dx, dy)$$

$$+ \int_E f(x) g(x) \kappa(dx),$$

where $\mathcal{E}^{(c)}$ is the strongly local part of $\mathcal{E}$.

For $u \in \mathcal{F}$, the martingale part $M_t^u$ in (1.2) can be decomposed as

$$M_t^u = M_t^{u,c} + M_t^{u,j} + M_t^{u,\kappa} \qquad \text{for every } t \in [0, \infty[,$$

$\mathbf{P}_x$-a.s. for q.e. $x \in E$, where $M_t^{u,c}$ is the continuous part of the martingale $M^u$ and

$$M_t^{u,j} = \lim_{\varepsilon \downarrow 0} \left\{ \sum_{0 < s \leq t} (u(X_s) - u(X_{s-})) \mathbf{1}_{\{|u(X_s) - u(X_{s-})| > \varepsilon\}} \mathbf{1}_{\{s < \zeta\}} \right.$$

$$-\int_0^t \left( \int_{\{y \in E : |u(y) - u(X_s)| > \varepsilon\}} (u(y) - u(X_s)) N(X_s, dy) \right) dH_s \bigg\},$$

$$M_t^{u,\kappa} = \int_0^t u(X_s) N(X_s, \{\Delta\}) \, dH_s - u(X_{\zeta-}) \mathbf{1}_{\{t \geq \zeta\}}$$

are the jump and killing parts of $M^u$, respectively. All three terms in this decomposition of $M^u$ are elements of $\overset{\circ}{\mathcal{M}}$; see Theorem A.3.9 of [7]. The limit in the expression for $M^{u,j}$ is in the sense of convergence in the norm of the space of MAF's of finite energy and of convergence in probability under $\mathbf{P}_x$ for q.e. $x \in E$ (see [7]).

Let $\mathcal{N}_c^* \subset \mathcal{N}_c$ denote the class of continuous additive functionals of the form $N^u + \int_0^\cdot g(X_s) \, ds$ for some $u \in \mathcal{F}$ and $g \in L^2(E; m)$. Nakao [14] constructed a linear map $\Gamma$ from $\overset{\circ}{\mathcal{M}}$ into $\mathcal{N}_c^*$ in the following way. It is shown in [14] that, for every $Z \in \overset{\circ}{\mathcal{M}}$, there is a unique $w \in \mathcal{F}$ such that

(1.4) $\qquad \mathcal{E}_1(w, f) = \frac{1}{2} \mu_{\langle M^f + M^{f,\kappa}, Z \rangle}(E) \qquad$ for every $f \in \mathcal{F}$.

This unique $w$ is denoted by $\gamma(Z)$. The operator $\Gamma$ is now defined by

(1.5) $\qquad \Gamma(Z)_t := N_t^{\gamma(Z)} - \int_0^t \gamma(Z)(X_s) \, ds \qquad$ for every $Z \in \overset{\circ}{\mathcal{M}}$.

Nakao showed that $\Gamma(Z)$ is characterized by the following equation

(1.6) $\qquad \lim_{t \downarrow 0} \frac{1}{t} \mathbf{E}_{g \cdot m}[\Gamma(Z)_t] = -\frac{1}{2} \mu_{\langle M^g + M^{g,\kappa}, Z \rangle}(E) \qquad$ for every $g \in \mathcal{F}_b$.

Here, $\mathcal{F}_b := \mathcal{F} \cap L^\infty(E; m)$. So, in particular, we have $\Gamma(M^u) = N^u$ for $u \in \mathcal{F}$. Nakao [14] then used the operator $\Gamma$ to define a stochastic integral

(1.7) $\qquad \int_0^t f(X_s) \, dN_s^u := \Gamma(f * M^u)_t - \frac{1}{2} \langle M^{f,c} + M^{f,j}, M^{u,c} + M^{u,j} \rangle_t,$

where $u \in \mathcal{F}$, $f \in \mathcal{F} \cap L^2(E; \mu_{\langle u \rangle})$ and $(f * M^u)_t := \int_0^t f(X_{s-}) \, dM_s^u$. If we define

$$\widetilde{\mathcal{N}}_c := \{N \in \mathcal{N}_c \mid N = N^u + A^\mu \text{ for some } u \in \mathcal{F}$$

and some signed smooth measure $\mu\}$,

then we see, by (1.5), that $\int_0^\cdot f(X_s) \, dN_s^u \in \widetilde{\mathcal{N}}_c$ if $u \in \mathcal{F}$ and $f \in \mathcal{F} \cap L^2(E; \mu_{\langle u \rangle})$. However, the conditions imposed on the integrand $f(X_t)$ and on the integrator $N^u$ in Nakao's stochastic integral are too restrictive for certain applications, in particular the perturbation theory of general symmetric Markov processes, which requires more general integrators as well as integrands; see [2].

The purpose of this paper is to provide a new way of defining $\Gamma(M)$ and Nakao's stochastic integral for zero-energy AF's $N^u$.

For a finite rcll AF $M_t$, it is known (see [3], Lemma 3.2) that there is a Borel function $\varphi$ on $E \times E$ with $\varphi(x,x) = 0$ for all $x \in E$ so that

$$(1.8) \qquad M_t - M_{t-} = \varphi(X_{t-}, X_t) \qquad \text{for every } t \in ]0, \zeta[, \mathbf{P}_m\text{-a.e.}$$

Such a $\varphi$ is uniquely determined up to $J$-negligible sets. We will call $\varphi$ the *jump function* of $M$. When $M = M^u$, $u \in \mathcal{F}$, the jump function $\varphi$ for $M^u$ can be taken to be as $\varphi(x,y) = u(y) - u(x)$ for $(x,y) \in E \times E$, with $u(\Delta) := 0$. We have a similar result for locally square-integrable MAF's on $[\![0, \zeta[\![$ [see Definition 2.5(iii) for the definition of a locally square-integrable MAF on $[\![0, \zeta[\![$]. Let $M$ be a locally square-integrable MAF on $[\![0, \zeta[\![$. There then exists a jump function $\varphi$ on $E \times E$ for $M$ satisfying the property (1.8) (see Corollary 2.9). Assume that

$$(1.9) \qquad \int_0^t \int_E (\widehat{\varphi}^2 \mathbf{1}_{\{|\widehat{\varphi}| \leq 1\}} + |\widehat{\varphi}| \mathbf{1}_{\{|\widehat{\varphi}| > 1\}})(X_s, y) N(X_s, dy) \, dH_s < \infty$$

for every $t < \zeta$,

$\mathbf{P}_x$-a.s. for q.e. $x \in E$, where $\widehat{\varphi}(x,y) := \varphi(x,y) + \varphi(y,x)$ for $x, y \in E$. By Lemma 3.2 below, there is a unique purely discontinuous local MAF $K$ on $[\![0, \zeta[\![$ with

$$K_t - K_{t-} = -\widehat{\varphi}(X_{t-}, X_t) \qquad \text{for } t < \zeta, \ \mathbf{P}_x\text{-a.s. for q.e. } x \in E.$$

Define $\mathbf{P}_m$-a.e. on $[0, \zeta[$,

$$\Lambda(M)_t := -\tfrac{1}{2}(M_t + M_t \circ r_t + \varphi(X_t, X_{t-}) + K_t) \qquad \text{for } t \in [0, \zeta[,$$

where $r_t$ is the time-reversal operator at time $t > 0$. Note that since $X$ is symmetric, the measure $\mathbf{P}_m$, when restricted to $\{t < \zeta\}$, is invariant under $r_t$. This time reversibility plays an important role in this paper. So, $\Lambda(M)$ is clearly well defined on $[\![0, \zeta[\![$ under the $\sigma$-finite measure $\mathbf{P}_m$. It will be shown in Theorem 2.18 and Remark 3.4(ii) below that $\Lambda(M)$ is a continuous even AF of $X$ on $[\![0, \zeta[\![$ admitting $m$-null set. Note that when $M = M^u$ for some $u \in \mathcal{F}$, $\varphi(x,y) = u(y) - u(x)$ is antisymmetric and so $\widehat{\varphi} = 0$. Thus, $\mathbf{P}_m$-a.e. on $\{t < \zeta\}$,

$$\Lambda(M^u)_t := -\tfrac{1}{2}(M_t^u + M_t^u \circ r_t + u(X_{t-}) - u(X_t)) = N_t^u.$$

The last identity follows by applying the time-reversal operator to both sides of (1.2) and using the fact that $N_t^u \circ r_t = N_t^u$ $\mathbf{P}_m$-a.e. on $[\![0, \zeta[\![$ (cf. [4], Theorem 2.1). It then follows for every $u \in \mathcal{F}$ that $\Lambda(M^u) = \Gamma(M^u)$ on $[\![0, \zeta[\![$ $\mathbf{P}_m$-a.e. We will show in Theorem 3.6 below that this holds when $M^u$ is replaced by any $M \in \overset{\circ}{\mathcal{M}}$. Therefore, under the $\sigma$-finite measure $\mathbf{P}_m$, $\Lambda$ is a genuine extension of Nakao's map $\Gamma$.

A function $f$ is said to be *locally in* $\mathcal{F}$ (denoted as $f \in \mathcal{F}_{\text{loc}}$) if there is an increasing sequence of finely open Borel sets $\{D_k, k \geq 1\}$ with $\bigcup_{k=1}^{\infty} D_k = E$ q.e. and for every $k \geq 1$, there is $f_k \in \mathcal{F}$ such that $f = f_k$ $m$-a.e. on $D_k$. For two subsets $A, B$ of $E$, we denote $A = B$ q.e. if $A \triangle B := (A \setminus B) \cup (B \setminus A)$ is $\mathcal{E}$-polar. By definition, every $f \in \mathcal{F}_{\text{loc}}$ admits a quasi-continuous $m$-version, so we may assume that all $f \in \mathcal{F}_{\text{loc}}$ are quasi-continuous. We then have $f = f_k$ q.e. on $D_k$. For $f \in \mathcal{F}_{\text{loc}}$, $M^{f,c}$ is well defined as a continuous MAF on $[\![0, \zeta[\![$ of locally finite energy. Moreover, for $f \in \mathcal{F}_{\text{loc}}$ and a locally square-integrable MAF $M$ on $[\![0, \zeta[\![$,

$$t \mapsto (f * M)_t := \int_0^t f(X_{s-})\, dM_s$$

is a locally square-integrable MAF on $[\![0, \zeta[\![$. Here, for a locally square-integrable MAF $M$ on $[\![0, \zeta[\![$, denote by $M^c$ its continuous part, which is also a locally square-integrable MAF on $[\![0, \zeta[\![$ (see Theorem 8.23 in [9]).

DEFINITION 1.2 (*Stochastic integral*). Suppose that $M$ is a locally square-integrable MAF on $[\![0, \zeta[\![$ and that $f \in \mathcal{F}_{\text{loc}}$. Let $\varphi: E \times E \to \mathbb{R}$ be a jump function for $M$ and assume that $\varphi$ satisfies condition (1.9). Define, on $[\![0, \zeta[\![$,

$$\int_0^t f(X_{s-})\, d\Lambda(M)_s$$
$$:= \Lambda(f * M)_t - \tfrac{1}{2} \langle M^{f,c}, M^c \rangle_t$$
$$+ \tfrac{1}{2} \int_0^t \int_E (f(y) - f(X_s))\varphi(y, X_s) N(X_s, dy)\, dH_s,$$

whenever $\Lambda(f * M)$ is well defined and the third term in the right-hand side of (3.10) is absolutely convergent.

The above stochastic integral is well defined on $[\![0, \zeta[\![$ under the $\sigma$-finite measure $\mathbf{P}_m$ and extends that of Nakao (1.7). [See Remark 3.9(i) and Theorem 3.10 below.] We will show in Theorem 4.7 below that it enjoys a generalized Itô formula.

**2. Additive functionals.** In this section, we will prove some facts about additive functionals, to be used later. We begin with some details on the completion of filtrations. Let $\mathcal{P}(E)$ be the family of all probability measures on $E$. For each $\nu \in \mathcal{P}(E)$, let $\mathcal{F}_\infty^\nu$ (resp., $\mathcal{F}_t^\nu$) be the $\mathbf{P}_\nu$-completion of $\mathcal{F}_\infty^0$ (resp., $\mathbf{P}_\nu$-completion of $\mathcal{F}_t^0$ in $\mathcal{F}_\infty^\nu$) and set $\mathcal{F}_\infty := \bigcap_{\nu \in \mathcal{P}(E)} \mathcal{F}_\infty^\nu$ and $\mathcal{F}_t := \bigcap_{\nu \in \mathcal{P}(E)} \mathcal{F}_t^\nu$. Let $\mathcal{F}_\infty^m$ (resp., $\mathcal{F}_t^m$) be the $\mathbf{P}_m$-completion of $\mathcal{F}_\infty^0$ (resp., $\mathbf{P}_m$-completion of $\mathcal{F}_t^0$ in $\mathcal{F}_\infty^m$). Although $m$ may not be a finite measure on $E$, we do have $\mathcal{F}_\infty \subset \mathcal{F}_\infty^m$, $\mathcal{F}_t \subset \mathcal{F}_t^m$ because for $g \in L^1(E; m)$ with $0 < g \leq 1$ on $E$ satisfying $gm \in \mathcal{P}(E)$, $\mathbf{P}_{gm}$-negligibility is the same as $\mathbf{P}_m$-negligibility.

For a fixed filtration $(\mathcal{M}_t)$ on $(\Omega, \mathcal{M})$, we recall the notions of $(\mathcal{M}_t)$-predictability, $(\mathcal{M}_t)$-optionality and $(\mathcal{M}_t)$-progressive measurability as follows (see [15] for more details). On $[0, \infty[ \times \Omega$, the $(\mathcal{M}_t)$-*predictable* [resp., $(\mathcal{M}_t)$-*optional*] $\sigma$-field $\mathcal{P}(\mathcal{M}_t)$ [resp., $\mathcal{O}(\mathcal{M}_t)$] is defined as the smallest $\sigma$-field over $[0, \infty[ \times \Omega$ containing all $\mathbf{P}_\nu(\mathcal{M})$-evanescent sets for all $\nu \in \mathcal{P}(E_\Delta)$ and with respect to which all $\mathcal{M}_t$-adapted lcrl (left-continuous, right-limited) (resp., rcll) processes are measurable. A process $\phi(s, \omega)$ on $[0, \infty[ \times \Omega$ is said to be $(\mathcal{M}_t)$-*progressively measurable* provided $[0, t] \times \Omega \ni (s, \omega) \mapsto \phi(s, \omega)$ is $\mathcal{B}([0, t]) \otimes \mathcal{M}_t$-measurable for all $t > 0$. It is well known that $(\mathcal{M}_t)$-predictability implies $(\mathcal{M}_t)$-optionality, which in turn implies $(\mathcal{M}_t)$-progressive measurability.

For $[0, \infty]$-valued functions $S, T$ on $\Omega$ with $S \leq T$, we employ the usual notation for *stochastic intervals*; for example,

$$[\![S, T[\![ := \{(t, \omega) \in [0, \infty[ \times \Omega \mid S(\omega) \leq t < T(\omega)\},$$

the other species of stochastic intervals being defined analogously. We write $[\![S]\!] := [\![S, S]\!]$ for the *graph* of $S$. Note that these are all subsets of $[0, \infty[ \times \Omega$. If $S$ and $T$ are $(\mathcal{M}_t)$-stopping times, then $[\![S, T]\!]$, $[\![S, T[\![, \ldots$ and $[\![S]\!]$ are $(\mathcal{M}_t)$-optional (see Theorem 3.16 in [9]).

DEFINITION 2.1 (AF). An $(\mathcal{F}_t)$-adapted [resp., $(\mathcal{F}_t^m)$-adapted] process $A = (A_t)_{t \geq 0}$ with values in $[-\infty, \infty]$ is said to be an *additive functional* (AF in short) (resp., *AF admitting $m$-null set*) if there exist a *defining set* $\Xi \in \mathcal{F}_\infty$ and an $\mathcal{E}$-polar (resp., $m$-null) set $N$ satisfying the following conditions:

(i) $\mathbf{P}_x(\Xi) = 1$ for all $x \in E \setminus N$;
(ii) $\theta_t \Xi \subset \Xi$ for all $t \geq 0$; in particular, $\omega_\Delta \in \Xi$ and $\mathbf{P}_\Delta(\Xi) = 1$ because $\omega_\Delta = \theta_{\zeta(\omega)}(\omega)$ for all $\omega \in \Xi$;
(iii) for all $\omega \in \Xi$, $A_\cdot(\omega)$ is right-continuous with left limits on $[0, \zeta(\omega)[$, $A_0(\omega) = 0, |A_t(\omega)| < \infty$ for $t < \zeta(\omega)$ and $A_{t+s}(\omega) = A_t(\omega) + A_s(\theta_t \omega)$ for all $t, s \geq 0$;
(iv) for all $t \geq 0$, $A_t(\omega_\Delta) = 0$; in particular, under the additivity in (iii), $A_t(\omega) = A_{\zeta(\omega)}(\omega)$ for all $t \geq \zeta(\omega)$ and $\omega \in \Xi$.

An AF $A$ (admitting $m$-null set) is called *right-continuous with left limits* (rcll AF in brief) if $A_{\zeta(\omega)-}$ exists for each $\omega \in \Xi$. An AF $A$ (admitting $m$-null set) is said to be *finite* [resp., *continuous additive functional* (CAF in brief)] if $|A_t(\omega)| < \infty$, $t \in [0, \infty[$ (resp. $t \mapsto A_t(\omega)$ is continuous on $[0, \infty[$) for each $\omega \in \Xi$. A $[0, \infty[$-valued CAF is called a *positive continuous additive functional* (PCAF in short). Two AF's $A$ and $B$ are called *equivalent* if there exists a common defining set $\Xi \in \mathcal{F}_\infty$ and an $\mathcal{E}$-polar set $N$ such that $A_t(\omega) = B_t(\omega)$ for all $t \in [0, \infty[$ and $\omega \in \Xi$. We call $A = (A_t)_{t \geq 0}$ an AF on $[\![0, \zeta[\![$ or a local AF (admitting $m$-null set) if $A$ is $(\mathcal{F}_t)$-adapted and satisfies (i), (ii), (iv) and the property (iii)′ in which (iii) is modified so that the additivity condition is required only for $t + s < \zeta(\omega)$. The notions of rcll AF, CAF and PCAF on $[\![0, \zeta[\![$ are defined similarly. Two AF's on $[\![0, \zeta[\![$, $A$ and $B$, are called *equivalent* if there exists a common defining set $\Xi \in \mathcal{F}_\infty$ and an $\mathcal{E}$-polar set $N$ such that $A_t(\omega) = B_t(\omega)$ for all $t \in [0, \zeta[$ and $\omega \in \Xi$.

REMARK 2.2. Any PCAF $A$ on $[\![0, \zeta[\![$ can be extended to a PCAF by setting

$$A_t(\omega) := \begin{cases} \lim_{u \uparrow \zeta} A_u(\omega), & \text{if } t \geq \zeta(\omega) > 0, \\ 0, & \text{if } t \geq \zeta(\omega) = 0, \end{cases}$$

for $\omega \in \Xi$ and setting $A_t(\omega) \equiv 0$ for $\omega \in \Xi^c$. The $(\mathcal{F}_t)$-adaptedness of this extended $A$ holds as follows: for a fixed $T > 0$, we know $\{A_t \leq T\} \cap \{t < \zeta\} \in \mathcal{F}_t$. From this, we have the $\mathcal{F}_\zeta$-measurability of $\{A_\zeta \leq T\}$. Indeed, $\{A_\zeta \leq T\} = \bigcap_{t \in \mathbb{Q}_+} \{A_t \leq T, t < \zeta\} \in \mathcal{F}_\zeta$ as $\{A_t \leq T, t < \zeta\} \in \mathcal{F}_\zeta$ for any $t \geq 0$. Thus, $\{A_t \leq T\} \cap \{t \geq \zeta\} = \{A_\zeta \leq T\} \cap \{t \geq \zeta\} \in \mathcal{F}_t$. Therefore, $\{A_t \leq T\} \in \mathcal{F}_t$ for any $T > 0$, which gives the $(\mathcal{F}_t)$-adaptedness of $A$. Noting that $\zeta \circ \theta_t = \zeta - t$ if $t < \zeta$ and $\zeta \circ \theta_t = 0$ if $t \geq \zeta$, we conclude that $A_\zeta = A_t + A_\zeta \circ \theta_t$ for any $t \in [0, \infty[$ on $\Xi$. Consequently, $A_{t+s} = A_t + A_s \circ \theta_t$ holds for any $t, s \in [0, \infty[$ on $\Xi$.

The following lemma is a special case of [14], Theorem 2.2.

LEMMA 2.3. *Let $A, B$ be PCAF's such that for $m$-a.e. $x \in E$, $\mathbf{E}_x[A_t] = \mathbf{E}_x[B_t]$ for all $t \geq 0$ and suppose that the Revuz measure $\mu_A$ has finite total mass. $A$ is then equivalent to $B$.*

REMARK 2.4. The above lemma may fail if the condition $\mu_A(E) < \infty$ is not satisfied. For example, take $E = \mathbb{R}^d$ with $d \geq 2$ and let $X$ be Brownian motion on $\mathbb{R}^d$ and $\mu_A(dx) = |x|^{-d-1} dx$. $\mu_A$ is then a smooth measure and corresponds to a PCAF $A$ of $X$. Let $B_t = A_t + t$, which is a PCAF of $X$ with Revuz measure $\mu_A(dx) + dx$. However,

$$\mathbf{E}_x[A_t] = \int_0^t \left( \int_{\mathbb{R}^d} p(s, x, y) |y|^{-d-1} dy \right) ds = \infty = \mathbf{E}_x[B_t]$$

for every $x \in \mathbb{R}^d \setminus \{0\}$.

Here, $p(s, x, y) = (2\pi t)^{-d/2} \exp(-|x - y|^2/(2t))$ is the transition density function of $X$.

As usual, if $T$ is an $(\mathcal{F}_t)$-stopping time and $M$ a process, then $M^T$ is the stopped process defined by $M_t^T := M_{t \wedge T}$. Following [9], we give the notion of local martingales of interval type.

DEFINITION 2.5 (*Processes of interval type*). Let $\mathcal{D}$ be a class of $(\mathcal{F}_t)$-adapted processes and denote by $\mathcal{D}_{\text{loc}}$ its localization (resp., by $\mathcal{D}_{f\text{-loc}}$ its localization by a nest of finely open Borel sets); that is, $M \in \mathcal{D}_{\text{loc}}$ (resp., $M \in \mathcal{D}_{f\text{-loc}}$) if and only if there exists a sequence $M^n \in \mathcal{D}$ and an increasing sequence of stopping times $T_n$ with $T_n \to \infty$ (resp., a nest $\{G_n\}$ of finely open Borel sets) such that $M^{T_n} = (M^n)^{T_n}$ (resp., $M_t = M_t^n$ for $t < \tau_{G_n}$) for each $n$. Here, a family $\{G_n\}$ of finely open Borel sets is called a *nest* if $\mathbf{P}_x(\lim_{n \to \infty} \tau_{G_n} = \zeta) = 1$ for q.e. $x \in E$. (However, see Lemma 3.1.) Clearly, $\mathcal{D} \subset \mathcal{D}_{\text{loc}}$ (resp., $\mathcal{D} \subset \mathcal{D}_{f\text{-loc}}$) and $(\mathcal{D}_{\text{loc}})_{\text{loc}} = \mathcal{D}_{\text{loc}}$ [resp., $(\mathcal{D}_{f\text{-loc}})_{f\text{-loc}} = \mathcal{D}_{f\text{-loc}}$]. If $\mathcal{D}$ is a subclass of AF's, then so is $\mathcal{D}_{\text{loc}}$ [for if $M \in \mathcal{D}_{\text{loc}}$, then there exist $M^n$ and $T_n$ as above and for each $\omega$ and $t, s \geq 0$, there exists $n \in \mathbb{N}$ with $s + t < T_n(\omega)$ and $s < T_n(\theta_t \omega)$, hence $M_{t+s}(\omega) = M_t(\omega) + M_s(\theta_t \omega)$], while $\mathcal{D}_{f\text{-loc}}$ is contained in the class of AF's on $[\![0, \zeta[\![$.

(i) $B \subset [0, \infty[ \times \Omega$ is called a *set of interval type* if there exists a nonnegative random variable $S$ such that for each $\omega \in \Omega$, the section $B_\omega := \{t \in [0, \infty[ \, | \, (t, \omega) \in B\}$ is $[0, S(\omega)]$ or $[0, S(\omega)[$ and $B_\omega \neq \varnothing$.

(ii) Let $B$ be an $(\mathcal{F}_t)$-optional set of interval type. A real-valued stochastic process $M$ on $B$ [i.e., $M \mathbf{1}_B = (M_t(\omega) \mathbf{1}_B(t, \omega))_{t \geq 0}$ is a real-valued stochastic process] is said to be *in $\mathcal{D}^B$* if and only if there exists $N \in \mathcal{D}$ such that $M \mathbf{1}_B = N \mathbf{1}_B$ and is said to be *locally in $\mathcal{D}$ on $B$* [write $M \in (\mathcal{D}_{\text{loc}})^B$] if and only if $S := D_{B^c}$ is the debut of $B^c$ and there exists an increasing sequence of $(\mathcal{F}_t)$-stopping times $\{S_n\}$ with $\lim_{n \to \infty} S_n = S$ and a sequence of $M^n \in \mathcal{D}$ such that $B_\omega \subset \bigcup_{n=1}^\infty [0, S_n(\omega)]$ $\mathbf{P}_x$-a.s. $\omega \in \Omega$ and $(M \mathbf{1}_B)^{S_n} = (M^n \mathbf{1}_B)^{S_n}$ for all $n \in \mathbb{N}$ and $t \geq 0$, $\mathbf{P}_x$-a.s. $\omega \in \Omega$ for q.e. $x \in E$. Clearly, $\mathcal{D}^B \subset (\mathcal{D}_{\text{loc}})^B$. Moreover, $\mathcal{D}^{B_2} \subset \mathcal{D}^{B_1}$ and $(\mathcal{D}_{\text{loc}})^{B_2} \subset (\mathcal{D}_{\text{loc}})^{B_1}$ for any pair of $(\mathcal{F}_t)$-optional sets $B_1, B_2$ of interval type with $B_1 \subset B_2$.

(iii) Let $B$ be an $(\mathcal{F}_t)$-optional set of interval type. We set

$$\mathcal{M}^1 := \{M \mid M \text{ is a finite rcll AF}, \mathbf{E}_x[|M_t|] < \infty,$$
$$\mathbf{E}_x[M_t] = 0 \text{ for } \mathcal{E}\text{-q.e. } x \in E \text{ and all } t \geq 0\}$$

and speak of an element of $(\mathcal{M}^1)^B$ [resp., $(\mathcal{M}^1_{\text{loc}})^B$] as being an *MAF on $B$* (resp., *a local MAF on $B$*). Similarly,

$$\mathcal{M} := \{M \mid M \text{ is a finite rcll AF}, \mathbf{E}_x[M_t^2] < \infty,$$
$$\mathbf{E}_x[M_t] = 0 \text{ for } \mathcal{E}\text{-q.e. } x \in E \text{ and all } t \geq 0\}$$

and an element of $\mathcal{M}^B$ [resp., $(\mathcal{M}_{\text{loc}})^B$] is a *square-integrable MAF on $B$* (resp., *locally square-integrable MAF on $B$*). We further set

$$\mathcal{M}^c := \{M \in \mathcal{M} \mid M \text{ is a CAF}\},$$
$$\mathcal{M}^d := \{M \in \mathcal{M} \mid M \text{ is a purely discontinuous AF}\}$$

and an element of $(\mathcal{M}_{\mathrm{loc}}^c)^B$ [resp., $(\mathcal{M}_{\mathrm{loc}}^d)^B$] is called a *locally square-integrable continuous MAF on $B$* (resp., *locally square-integrable purely discontinuous MAF on $B$*). For $M \in (\mathcal{M}_{\mathrm{loc}})^B$, $M$ admits a unique decomposition $M = M^c + M^d$ with $M^c \in (\mathcal{M}_{\mathrm{loc}}^c)^B$ and $M^d \in (\mathcal{M}_{\mathrm{loc}}^d)^B$ (see Theorem 8.23 in [9]). In these definitions, we omit the usage "on $B$" when $B = [0, \infty[ \times \Omega$.

For a $[0, \infty]$-valued function $R$ on $\Omega$ and $A \subset \Omega$, $R_A := R \cdot \mathbf{1}_A + (+\infty) \cdot \mathbf{1}_{A^c}$ is called the *restriction* of $R$ on $A$. Clearly, $R \leq R_A$.

REMARK 2.6.  When $B = [\![0, R[\![$ for a given $(\mathcal{F}_t)$-stopping time $R$, there is another notion of "locally in $\mathcal{D}$ on $B$," obtained by replacing $(M\mathbf{1}_B)^{S_n} = (M^n\mathbf{1}_B)^{S_n}$ with $M^{S_n}\mathbf{1}_B = (M^n)^{S_n}\mathbf{1}_B$ in our definition; this is a weaker notion than ours because $t \mapsto \mathbf{1}_B(t, \omega)$ is decreasing and $\mathbf{1}_B(t, \omega)\mathbf{1}_B(s, \omega) = \mathbf{1}_B(t, \omega)$ for $s \leq t$ and $\omega \in \Omega$. This weaker notion is described in [15].

DEFINITION 2.7 (*MAF locally of finite energy*).  Recall that $\overset{\circ}{\mathcal{M}}$ is the totality of MAF's of finite energy, that is,

$$\overset{\circ}{\mathcal{M}} := \left\{ M \in \mathcal{M} \,\Big|\, e(M) := \lim_{t \downarrow 0} \frac{1}{2t} \mathbf{E}_m[M_t^2] < \infty \right\}.$$

We say that an AF $M$ on $[\![0, \zeta[\![$ is *locally in* $\overset{\circ}{\mathcal{M}}$ (and write $M \in \overset{\circ}{\mathcal{M}}_{f\text{-loc}}$) if there exists a sequence $\{M^n\}$ in $\overset{\circ}{\mathcal{M}}$ and a nest $\{G_n\}$ of finely open Borel sets such that $M_t = M_t^n$ for $t < \tau_{G_n}$ for each $n \in \mathbb{N}$. In case $X$ is a diffusion process with no killing inside $E$, we can define the predictable quadratic variation $\langle M \rangle$ for $M \in \overset{\circ}{\mathcal{M}}_{f\text{-loc}}$ as follows. First, note that $M_{t \wedge \tau_{G_n}}^n = M_{t \wedge \tau_{G_n}}^m$ for $n < m$ because of the continuity of $M^n$. Owing to the uniqueness of Doob–Meyer decomposition, we see that $\langle M^n \rangle_{t \wedge \tau_{G_n}} = \langle M^m \rangle_{t \wedge \tau_{G_n}}$. The predictable quadratic variation $\langle M \rangle$ of $M \in \overset{\circ}{\mathcal{M}}_{f\text{-loc}}$ as a PCAF is well defined by setting $\langle M \rangle_t = \langle M^n \rangle_t, t < \tau_{G_n}, n \in \mathbb{N}$, with Remark 2.2 and by choosing an appropriate defining set and $\mathcal{E}$-polar set of $\langle M \rangle$, where $M^n \in \overset{\circ}{\mathcal{M}}$ and $\{G_n\}$ is a nest of finely open Borel sets such that $M_t = M_t^n, t < \tau_{G_n}$.

PROPOSITION 2.8.  $(\mathcal{M}_{\mathrm{loc}})^{[\![0, \zeta[\![} \subset \overset{\circ}{\mathcal{M}}_{f\text{-loc}}$. *More precisely, for each* $M \in (\mathcal{M}_{\mathrm{loc}})^{[\![0, \zeta[\![}$, *there exists a nest* $\{G_k\}$ *of finely open Borel sets such that* $\mathbf{1}_{G_k} * M \in \overset{\circ}{\mathcal{M}}$ *for each* $k \in \mathbb{N}$ *and the predictable quadratic variation process* $\langle M \rangle$ *can be constructed as a PCAF.*

PROOF.  Let $M \in (\mathcal{M}_{\mathrm{loc}})^{[\![0, \zeta[\![}$. There then exists an increasing sequence $\{T_n\}$ of stopping times with $\lim_{n \to \infty} T_n = \zeta$ ($\mathbf{P}_x$-a.s. $\omega \in \Omega$ for q.e. $x \in E$) and $M^n \in \mathcal{M}_{\mathrm{loc}}$ such that $M_{t \wedge T_n} \mathbf{1}_{[0, \zeta[}(t \wedge T_n) = M_{t \wedge T_n}^n \mathbf{1}_{[0, \zeta[}(t \wedge T_n)$ holds for all $t \geq 0$ $\mathbf{P}_x$-a.s. for q.e. $x \in E$. We may assume that it holds for all $\omega \in \Omega$ by changing the sample space. Note that $[0, \zeta(\omega)[ \subset \bigcup_{n=1}^{\infty} [0, T_n(\omega)]$ for all $\omega \in \Omega$. Hence, $M^m_{t \wedge T_n} \mathbf{1}_{[0,\zeta[}(t \wedge T_n) = M^n_{t \wedge T_n} \mathbf{1}_{[0,\zeta[}(t \wedge T_n)$ for $n < m$. As noted in Definition 2.5, we see that $M$ is an AF on $[\![0, \zeta[\![$. Owing to the uniqueness of the Doob–Meyer decomposition for semimartingales on $[\![0, \zeta[\![$ (see [9]), we have $\langle M^m \rangle_{t \wedge T_n} \mathbf{1}_{[0,\zeta[}(t \wedge T_n) = \langle M^n \rangle_{t \wedge T_n} \mathbf{1}_{[0,\zeta[}(t \wedge T_n)$ for $n < m$. Thus, we have $\langle M^m \rangle_t = \langle M^n \rangle_t$ for $t < T_n$ and $n < m$. The predictable quadratic variation $\langle M \rangle$ of $M$ is therefore well defined by setting $\langle M \rangle_t := \langle M^n \rangle_t$ for $t < T_n$. Setting $\langle M \rangle_t := \langle M \rangle_\zeta := \lim_{s \uparrow \zeta} \langle M \rangle_s$ for all $t \geq \zeta$, we obtain a PCAF because of Remark 2.2. Let $\mu_{\langle M \rangle}$ be the Revuz measure corresponding to $\langle M \rangle$ and $\{F_k\}$ an $\mathcal{E}$-nest of closed sets such that $\mu_{\langle M \rangle}(F_k) < \infty$ for each $k$, and let $G_k$ be the fine interior of $F_k$. $\{G_k\}$ is then a nest. In view of the proofs of Theorem 5.6.1 and Lemma 5.6.2 in [7], the stochastic integral $\mathbf{1}_{G_k} * M$ is of finite energy with $\mathbf{e}(\mathbf{1}_{G_k} * M) = \frac{1}{2} \mu_{\langle M \rangle}(G_k)$ and its predictable quadratic variation $\langle \mathbf{1}_{G_k} * M \rangle$ is a PCAF. Let $\mu_k$ (resp., $\overset{\circ}{\mu}_k$) be the Revuz measure corresponding to $\langle \mathbf{1}_{G_k} * M \rangle$ (resp., $\langle \mathbf{1}_{G_k} * M, M \rangle$). By Lemma 5.6.2 in [7], for $M_i \in \overset{\circ}{\mathcal{M}}$ and $f_i \in L^2(E; d\mu_{\langle M_i \rangle})$ $(i = 1, 2)$, we have $f_1 f_2 d\mu_{\langle M_1, M_2 \rangle} = d\mu_{\langle f_1 * M_1, f_2 * M_2 \rangle}$, hence $\int_0^t (f_1 f_2)(X_s) \, d\langle M_1, M_2 \rangle_s = \langle f_1 * M_1, f_2 * M_2 \rangle_t$. From this, we see that $\langle \mu_k, f^2 \rangle = \langle \overset{\circ}{\mu}_k, f^2 \rangle = \langle \mathbf{1}_{G_k} \mu_{\langle M \rangle}, f^2 \rangle$ for any $f \in L^2(E; \mu_{\langle M \rangle})$; consequently, we have $\mu_k = \overset{\circ}{\mu}_k = \mathbf{1}_{G_k} \mu_{\langle M \rangle}$ by $\mu_{\langle M \rangle}(G_k) < \infty$. This yields $\langle \mathbf{1}_{G_k} * M \rangle_t = \langle \mathbf{1}_{G_k} * M, M \rangle_t = \int_0^t \mathbf{1}_{G_k}(X_s) \, d\langle M \rangle_s$ for $t < \zeta$, hence $\langle M - \mathbf{1}_{G_k} * M \rangle_t = 0$ for $t < \tau_{G_k}$. Therefore, $M_t = (\mathbf{1}_{G_k} * M)_t$ for $t < \tau_{G_k}$ and $\mathbf{1}_{G_k} * M \in \overset{\circ}{\mathcal{M}}$. $\square$

COROLLARY 2.9. *Let $M$ be a locally square-integrable MAF on $[\![0, \zeta[\![$, that is, $M \in (\mathcal{M}_{\text{loc}})^{[\![0,\zeta[\![}$. There then exists a Borel function $\varphi$ on $E \times E$ with $\varphi(x, x) = 0$ for all $x \in E$ such that*

$$M_t - M_{t-} = \varphi(X_{t-}, X_t) \qquad \text{for every } t \in ]0, \zeta[, \; \mathbf{P}_m\text{-a.e.}$$

PROOF. By the proof of Proposition 2.8, there exists an $\mathcal{E}$-nest $\{F_k\}$ such that for each $k \in \mathbb{N}$ $M^k := \mathbf{1}_{F_k} * M \in \overset{\circ}{\mathcal{M}}$ and $M_t = M^k_t$, $t < \tau_{F_k}$. Let $\varphi_k$ be the jump function corresponding to $M^k$. We then have $\varphi_k(X_{t-}, X_t) = \varphi_\ell(X_{t-}, X_t)$, $t < \tau_{F_k}$, $\mathbf{P}_m$-a.e., for $k < \ell$. From this, we see that $\varphi_k = \varphi_\ell$ $J$-a.e. on $F_k \times F_k$. We construct a Borel function $\varphi$ on $E \times E$ in the following manner. We set $F_0 := \varnothing$, $\varphi(x, y) := \varphi_k(x, y)$ for $(x, y) \in F_k \times F_k \setminus (F_{k-1} \times F_{k-1})$, $k \in \mathbb{N}$, $\varphi(x, y) := 0$ if $(x, y) \in E \times E \setminus (\bigcup_{k=1}^\infty F_k \times \bigcup_{k=1}^\infty F_k)$. $\varphi$ then satisfies $\varphi(x, x) = 0$ for $x \in E$. We also have $\varphi = \varphi_k$ $J$-a.e. on $F_k \times F_k$. Consequently, $\varphi(X_{t-}, X_t) = \varphi_k(X_{t-}, X_t)$, $t < \tau_{F_k}$, $\mathbf{P}_m$-a.e. This means that $M_t - M_{t-} = \varphi(X_{t-}, X_t)$, $t < \tau_{F_k}$, $\mathbf{P}_m$-a.e. Therefore, $M_t - M_{t-} = \varphi(X_{t-}, X_t)$, $0 < t < \zeta$ $\mathbf{P}_m$-a.e. $\square$

We recall the definition of the shift operator $\theta_s$ and the time-reversal operator $r_t$ on the path space $\Omega$. For each $s \geq 0$, the shift operator $\theta_s$ is defined by $\theta_s \omega(t) := \omega(t+s)$ for $t \in [0, \infty[$. Given a path $\omega \in \{t < \zeta\}$, the operator $r_t$ is defined by

$$(2.1) \qquad r_t(\omega)(s) := \begin{cases} \omega((t-s)-), & \text{if } 0 \leq s \leq t, \\ \omega(0), & \text{if } s \geq t. \end{cases}$$

Here, for $r > 0$, $\omega(r-) := \lim_{s \uparrow r} \omega(s)$ is the left limit at $r$ and we use the convention that $\omega(0-) := \omega(0)$. For a path $\omega \in \{t \geq \zeta\}$, we set $r_t(\omega) := \omega_\Delta$. We note that

$$(2.2) \qquad \begin{aligned} \lim_{s \downarrow 0} r_t(\omega)(s) &= \omega(t-) = r_t(\omega)(0) \quad \text{and} \\ \lim_{s \uparrow t} r_t(\omega)(s) &= \omega(0) = r_t(\omega)(t). \end{aligned}$$

A key consequence of the $m$-symmetry assumption on the Hunt process $X$ is that the measure $\mathbf{P}_m$, when restricted to $\{t < \zeta\}$, is invariant under the time-reversal operator $r_t$.

Clearly for $t, s > 0$, $\theta_s : \Omega \to \Omega$ is $\mathcal{F}_{t+s}^m/\mathcal{F}_t^m$-measurable. The following lemma deals with the measurability issue of the time-reversal operator $r_t$.

LEMMA 2.10. *For each $t > 0$, $r_t : \Omega \to \Omega$ is $\mathcal{F}_t^0/\mathcal{F}_\infty^0$-measurable and $\mathcal{F}_t^m/\mathcal{F}_t^m$-measurable.*

PROOF. Let $F_i \in \mathcal{B}(E_\Delta)$ and $s_i \in [0, \infty[$, $i = 1, 2, \ldots, n$, with $s_1 < s_2 < \cdots < s_k \leq t < s_{k+1} < \cdots < s_n$ for some $k \in \{1, 2, \ldots, n\}$. Then, $r_t^{-1}(\bigcap_{i=1}^n X_{s_i}^{-1}(F_i)) = \bigcap_{i=1}^n (X_{s_i} \circ r_t)^{-1}(F_i)$ is equal to $\bigcap_{i=1}^k (\{X_{(t-s_i)-} \in F_i, t < \zeta\} \cup \{\Delta \in F_i, t \geq \zeta\}) \cap \bigcap_{i=k+1}^n (\{X_0 \in F_i, t < \zeta\} \cup \{\Delta \in F_i, t \geq \zeta\}) \in \mathcal{F}_t^0$. Next, we show the $\mathcal{F}_t^m/\mathcal{F}_t^m$-measurability of $r_t$. Take $C \in \mathcal{F}_t^m$. There then exist $D \in \mathcal{F}_t^0$ and $N \in \mathcal{F}_\infty^0$ such that $C \triangle D \subset N$ and $\mathbf{P}_m(N) = 0$. Since $\mathbf{P}_m(\{\omega_\Delta\}) = 0$, by deleting $\{\omega_\Delta\} = \{\omega \in \Omega \mid \zeta(\omega) = 0\} \in \mathcal{F}_0^0 \subset \mathcal{F}_t^0$, we may assume that $\omega_\Delta \notin C \cup D \cup N$. Then, $r_t^{-1}(C) \triangle r_t^{-1}(D) \subset r_t^{-1}(N)$, $r_t^{-1}(D)$, $r_t^{-1}(N) \in \mathcal{F}_t^0$ and $\mathbf{P}_m(r_t^{-1}(N)) = \mathbf{P}_m(r_t^{-1}(N) \cap \{t < \zeta\}) + \mathbf{1}_N(\omega_\Delta)\mathbf{P}_m(t \geq \zeta) = \mathbf{P}_m(N \cap \{t < \zeta\}) = 0$. □

DEFINITION 2.11. For any $t > 0$, we say that two sample paths $\omega$ and $\omega'$ are *t-equivalent* if $\omega(s) = \omega'(s)$ for all $s \in [0, t]$. We say that two sample paths $\omega$ and $\omega'$ are *pre-t-equivalent* if $\omega(s) = \omega'(s)$ for all $s \in [0, t[$.

For an rcll AF $A_t$ adapted to $(\mathcal{F}_t^0)_{t \geq 0}$, $A_t(\omega) = A_t(\omega')$ if $\omega$ and $\omega'$ are $t$-equivalent $A_{t-}(\omega) = A_{t-}(\omega')$ if $\omega$ and $\omega'$ are pre-$t$-equivalent. These conclusions may fail to hold if the measurability conditions are not satisfied. We need the following notion.

DEFINITION 2.12 (PrAF). A process $A = (A_t)_{t \geq 0}$ with values in $\overline{\mathbb{R}} := [-\infty, \infty]$ is said to be a *progressively additive functional* (PrAF in short) (resp., *PrAF admitting m-null set*) if $A$ is $(\mathcal{F}_t)$-adapted [resp., $(\mathcal{F}_t^m)$-adapted] and there exist *defining sets* $\Xi \in \mathcal{F}_\infty$, $\Xi_t \in \mathcal{F}_t$ [resp., $\Xi \in \mathcal{F}_\infty^m$, $\Xi_t \in \mathcal{F}_t^m$] for each $t > 0$ and an $\mathcal{E}$-polar (resp., $m$-null) set $N$ satisfying the following conditions:

(i) $\mathbf{P}_x(\Xi) = 1$ for all $x \in E \setminus N$, $\Xi \subset \Xi_t \subset \Xi_s$ for every $t > s > 0$ and $\Xi = \bigcap_{t>0} \Xi_t$;

(ii) $\theta_t \Xi \subset \Xi$ for all $t \geq 0$ and $\theta_{t-s}(\Xi_t) \subset \Xi_s$ for all $s \in ]0, t[$, and, in particular, $\omega_\Delta \in \Xi \subset \Xi_t$ and $\mathbf{P}_\Delta(\Xi) = \mathbf{P}_\Delta(\Xi_t) = 1$ under (i);

(iii) for all $\omega \in \Xi_t$, $A(\omega)$ is defined on $[0, t[$, is right continuous on $[0, t \wedge \zeta(\omega)[$ and has left limit on $]0, t] \cap ]0, \zeta(\omega)[$ such that $A_0(\omega) = 0$, $|A_s(\omega)| < \infty$ for $s \in [0, t \wedge \zeta(\omega)[$ and $A_{p+q}(\omega) = A_p(\omega) + A_q(\theta_p \omega)$ for all $p, q \geq 0$ with $p + q < t$;

(iv) for all $t \geq 0$, $A_t(\omega_\Delta) = 0$;

(v) for any $t > 0$ and pre-$t$-equivalent paths $\omega, \omega' \in \Omega$, $\omega \in \Xi_t$ implies that $\omega' \in \Xi_t$, $A_s(\omega) = A_s(\omega')$ for any $s \in [0, t[$ and $A_{s-}(\omega) = A_{s-}(\omega')$ for any $s \in ]0, t]$.

Furthermore, $A$ is called an *rcll PrAF* (or an *rcll PrAF admitting m-null set*) if, for each $t > 0$ and $\omega \in \Xi_t$, $s \mapsto A_s(\omega)$ is right-continuous on $[0, t[$ and has left-hand limits on $]0, t]$ and a PrAF (or a PrAF admitting $m$-null set) is said to be *finite* (resp., *continuous*) if $|A_s(\omega)| < \infty$ for all $s \in [0, t[$ (resp., continuous on $[0, t[$) for every $\omega \in \Xi_t$.

We say that an AF $A$ on $[\![0, \zeta[\![$ (resp., AF $A$ on $[\![0, \zeta[\![$ admitting $m$-null set) is a *PrAF on* $[\![0, \zeta[\![$ (resp., *PrAF on* $[\![0, \zeta[\![$ *admitting m-null set*) if $A$ is $(\mathcal{F}_t)$-adapted [resp., $(\mathcal{F}_t^m)$-adapted] and there exist $\Xi \in \mathcal{F}_\infty$, $\Xi_t \in \mathcal{F}_t$ (resp., $\Xi \in \mathcal{F}_\infty^m$, $\Xi_t \in \mathcal{F}_t^m$) for each $t > 0$ and an $\mathcal{E}$-polar (resp., $m$-null) set $N$ such that (i'), (ii), (iii'), (iv) and (v') hold—(i'): $\mathbf{P}_x(\Xi) = 1$ for all $x \in E \setminus N$, $\Xi \subset \Xi_t$ for all $t > 0$, $\Xi = \bigcap_{t>0} \Xi_t$ and $\Xi_t \cap \{t < \zeta\} \subset \Xi_s \cap \{s < \zeta\}$ for $s < t$; (iii'): for each $\omega \in \Xi_t \cap \{t < \zeta\}$, the same conclusion as in (iii) holds; (v'): for any $t > 0$ and pre-$t$-equivalent paths $\omega, \omega' \in \Omega \cap \{t < \zeta\}$, the same conclusion as in (v) holds.

The notion of *rcll PrAF on* $[\![0, \zeta[\![$ (*or rcll PrAF admitting m-null set*) is similarly defined.

REMARK 2.13.   (i) Our notion of PrAF is different from what is found in Walsh [16].

(ii) Every PrAF (resp., PrAF on $[\![0, \zeta[\![$) is an AF (resp., AF on $[\![0, \zeta[\![$).

(iii) The MAF $M^u$ and the CAF $N^u$ of 0-energy appearing in Fukushima's decomposition (1.2) can be regarded as finite rcll PrAF's in view of the proof of Theorem 5.2.2 in [7]. In this case, the defining sets for $M^u$ as PrAF are given by

$$\Xi := \{\omega \in \Omega \mid M^{u_n}_{s-}(\omega) \text{ converges uniformly on } ]0,t] \text{ for } \forall t \geq 0$$
$$\text{for some subsequence } n_k\} \in \mathcal{F}_\infty,$$
$$\Xi_t := \{\omega \in \Omega \mid M^{u_n}_{s-}(\omega) \text{ converges uniformly on } ]0,t]$$
$$\text{for some subsequence } n_k\} \in \mathcal{F}_t$$

for every $t > 0$, where $M^{u_n}_t := u_n(X_t) - u_n(X_0) - \int_0^t (u_n(X_s) - f_n(X_s))\,ds$ with $f_n := n(u - nR_{n+1}u)$ and $u_n := R_1 f_n = nR_{n+1}u$. Hence, an MAF of stochastic integral type $\int_0^t g(X_{s-})\,dM^u_s$ $[g, u \in \mathcal{F}$ with $g \in L^2(E; \mu_{\langle u \rangle})]$ can be regarded as a finite rcll PrAF. Consequently, any MAF of finite energy can also be regarded as an rcll PrAF, in view of the assertion of Lemma 5.6.3 in [7] and Lemma 2.14 below.

(iv) Every $M \in \overset{\circ}{\mathcal{M}}_{f\text{-loc}}$ can be regarded as a PrAF on $[\![0,\zeta[\![$, hence every $M \in \mathcal{M}_{\text{loc}}^{[\![0,\zeta[\![}$ is also. Since every local martingale can be written as the sum of a local martingale with bounded jumps (and hence a locally square-integrable martingale) and a local martingale of finite variation, we conclude that every local MAF is a PrAF.

LEMMA 2.14. *Let $(A^n)$ be a sequence of finite rcll PrAF's with defining sets $\Xi^n \in \mathcal{F}_\infty$ and $\Xi^n_t \in \mathcal{F}_t$. For each $t > 0$, set*

$$\Xi_t := \left\{\omega \in \bigcap_{n \in \mathbb{N}} \Xi^n_t \,\Big|\, A^n \text{ converges uniformly on } [0,t[ \right\} \in \mathcal{F}_t$$

*and*

$$\Xi := \left\{\omega \in \bigcap_{n \in \mathbb{N}} \Xi^n \,\Big|\, A^n \text{ converges uniformly on } [0,t[ \text{ for every } t \in [0,\infty[ \right\} \in \mathcal{F}_\infty.$$

*Suppose that there exists an $\mathcal{E}$-polar set $N$ such that $\mathbf{P}_x(\Xi) = 1$ for $x \in E \setminus N$. If we define $A_t := \varliminf_{n \to \infty} A^n_t$ on $\Omega$, then $A$ is a finite rcll PrAF with its defining sets $\Xi$, $\Xi_t$.*

PROOF. We only show that for any $t > 0$ and pre-$t$-equivalent paths $\omega, \omega'$, $\omega \in \Xi_t$ implies that $\omega' \in \Xi_t$. Suppose that $\omega \in \Xi_t$ and $\omega$ is pre-$t$-equivalent to $\omega'$. It easy to see that $\omega' \in \bigcap_{n \in \mathbb{N}} \Xi^n_t$. We then see the uniform convergence of $A^n_{s-}(\omega') = A^n_{s-}(\omega)$ for $s \in ]0,t]$. Therefore, $\omega' \in \Xi_t$, $A_s(\omega') = A_s(\omega)$ for $s \in [0,t[$ and $A_{s-}(\omega') = A_{s-}(\omega)$ for $s \in ]0,t]$. □

Recall that $\{\theta_t, t > 0\}$ denotes the time-shift operators on the path space for the process $X$.

LEMMA 2.15.  *For $t, s > 0$:*

(i) $\theta_t r_{t+s}\omega$ *is $s$-equivalent to $r_s\omega$ if $t + s < \zeta(\omega)$ or $s \geq \zeta(\omega)$;*

(ii) $r_t \theta_s \omega$ *is pre-$t$-equivalent to $r_{t+s}\omega$ and, moreover, if $\omega$ is continuous at $s$, then $r_t\theta_s\omega$ is $t$-equivalent to $r_{t+s}\omega$.*

PROOF.  (i) We may assume that $t + s < \zeta(\omega)$. For $v \in [0, s]$,
$$\theta_t r_{t+s}\omega(v) = \omega((s-v)-) = r_s\omega(v)$$
and so $\theta_t r_{t+s}\omega$ is $s$-equivalent to $r_s\omega$.

(ii) Note that $t + s < \zeta(\omega)$ is equivalent to $t < \zeta(\theta_s\omega)$. It follows from the definition, if $t + s < \zeta(\omega)$, that

$$(2.3) \qquad (r_t\theta_s\omega)(v) = \begin{cases} \omega((t+s-v)-), & \text{if } 0 \leq v < t, \\ \omega(s), & \text{if } v = t, \end{cases}$$

while $r_{t+s}\omega(v) = \omega((t+s-v)-)$ for $0 \leq v \leq t$. Hence, typically, $r_t\theta_s\omega$ is only pre-$t$-equivalent to $r_{t+s}\omega$.  $\square$

Fix $t > 0$. Set $\mathcal{H}^t_s := \mathcal{F}_t$ for $s \in [0, t]$ and $\mathcal{H}^t_s := \mathcal{F}_s$ for $s \in ]t, \infty[$. $(\mathcal{H}^t_s)_{s \geq 0}$ is then a filtration over $(\Omega, \mathcal{F}_\infty)$ and $\mathcal{F}_s \subset \mathcal{H}^t_s$ for all $s \geq 0$.

LEMMA 2.16.  *The following assertions hold for any fixed $t > 0$:*

(i) *if we let $\varphi$ be a Borel function on $E \times E$ and set $X_{0-} := X_0$, then $[0, \infty[ \times \Omega \ni (s, \omega) \mapsto \mathbf{1}_{[\![0, \zeta[\![}(s, \omega) \mathbf{1}_{\Gamma_t}(\omega) \varphi(X_{s-}(\omega), X_s(\omega))$ is $(\mathcal{H}^t_s)$-optional for any $\Gamma_t \in \mathcal{F}_t$;*

(ii) *if we let $A$ be an rcll PrAF with defining sets $\Xi \in \mathcal{F}_\infty$, $\Xi_t \in \mathcal{F}_t$ and we set $A_{0-}(\omega) := 0$ and $A^t_s(\omega) := \mathbf{1}_{\Xi_t}(\omega)(\mathbf{1}_{[0,t]}(s)A_s(\omega) + \mathbf{1}_{]t,\infty[}(s)A_t(\omega))$ for $\omega \in \Omega$, then $[0, \infty[ \times \Omega \ni (s, \omega) \mapsto \mathbf{1}_{[\![0,\zeta[\![}(s,\omega)(A^t_s(\omega) - A^t_{s-}(\omega))$ is $(\mathcal{H}^t_s)$-optional.*

PROOF.  (i) Note that $\mathbf{1}_{[\![0,\zeta[\![}$ is $(\mathcal{H}^t_s)$-predictable. The assertion is clear if $\varphi = f \otimes g$ for bounded Borel functions $f, g$ on $E_\Delta$. The monotone class theorem for functions gives us the desired result.

(ii) Since $A^t$ is $(\mathcal{H}^t_s)$-adapted and rcll on $\Omega$ and $A^t_-$ is $(\mathcal{H}^t_s)$-adapted and lcrl on $\Omega$, $(s, \omega) \mapsto A_s(\omega)$ is $(\mathcal{H}^t_s)$-optional and $(s, \omega) \mapsto A^t_{s-}(\omega)$ is $(\mathcal{H}^t_s)$-predictable. Consequently, $(s, \omega) \mapsto A^t_s(\omega) - A^t_{s-}(\omega)$ is $(\mathcal{H}^t_s)$-optional.  $\square$

By Lemma 3.2 of [3], for a finite rcll AF $A = (A_t)_{t \geq 0}$, there is a Borel function $\varphi : E \times E \to \mathbb{R}$ with $\varphi(x, x) = 0$ for all $x \in E$ such that

$$(2.4) \qquad A_t - A_{t-} = \varphi(X_{t-}, X_t) \qquad \text{for every } t \in ]0, \zeta[, \mathbf{P}_m\text{-a.e.}$$

Moreover, if $\tilde{\varphi}$ is another such function, then $J(\varphi \neq \tilde{\varphi}) = 0$. As before, we refer to such a function $\varphi$ as a jump function for $A$. Recall that if $M \in \mathcal{M}^{[\![0,\zeta[\![}_{\text{loc}}$, then there exists a jump function $\varphi$ (unique in the above sense) so that $M_t - M_{t-} = \varphi(X_{t-}, X_t)$ for $t \in ]0, \zeta[, \mathbf{P}_m$-a.e.

LEMMA 2.17. *Let $A$ be a finite rcll PrAF with defining sets $\{\Xi, \Xi_t, t \geq 0\}$. There then exists a real-valued Borel function $\varphi$ on $E_\Delta \times E_\Delta$ with $\varphi(x, x) = 0$ for $x \in E_\Delta$ such that $A$ with defining sets*

$$\widetilde{\Xi} := \{\omega \in \Xi \mid A_s(\omega) - A_{s-}(\omega) = \varphi(X_{s-}(\omega), X_s(\omega)) \text{ for } s \in {]}0, \zeta(\omega){[}\},$$

$$\widetilde{\Xi}_t := \{\omega \in \Xi_t \mid A_s(\omega) - A_{s-}(\omega) = \varphi(X_{s-}(\omega), X_s(\omega)) \text{ for } s \in {]}0, t{[} \cap {]}0, \zeta(\omega){[}\}$$

*is again an rcll PrAF admitting m-null set. The analogous assertion holds for PrAF's on $[\![0, \zeta[\![$ and, in particular, for elements of $(\mathcal{M}_{\mathrm{loc}})^{[\![0,\zeta[\![}$.*

PROOF. Let $\varphi : E_\Delta \times E_\Delta \to \mathbb{R}$ be a Borel function vanishing on the diagonal and define $\widetilde{\Xi}$, $\widetilde{\Xi}_t$ in terms of $\varphi$, as above. Clearly, $\widetilde{\Xi} = \bigcap_{t>0} \widetilde{\Xi}_t$, $\widetilde{\Xi}_t \subset \widetilde{\Xi}_s$ for $s < t$. Moreover, we see that $\theta_t \widetilde{\Xi} \subset \widetilde{\Xi}$ for $t \geq 0$ and $\theta_{t-s}(\widetilde{\Xi}_t) \subset \widetilde{\Xi}_s$ for $s < t$. For two pre-$t$-equivalent paths $\omega, \omega'$, we see that $\omega \in \widetilde{\Xi}_t$ implies that $\omega' \in \widetilde{\Xi}_t$.

By the previous lemma,

$$\Gamma := \{(s, \omega) \mid \mathbf{1}_{[\![0, \zeta[\![}(s, \omega) \mathbf{1}_{\Xi_t}(\omega)(A_s^t(\omega) - A_{s-}^t(\omega) - \varphi(X_{s-}(\omega), X_s(\omega))) \neq 0\}$$

is $(\mathcal{H}_s^t)$-progressively measurable for any fixed $t > 0$ and the debut of $\Gamma$ is

$$D_\Gamma(\omega) := \inf\{s \geq 0 \mid \mathbf{1}_{[\![0, \zeta[\![}(s, \omega) \mathbf{1}_{\Xi_t}(\omega)(A_s^t(\omega) - A_{s-}^t(\omega)$$
$$- \varphi(X_{s-}(\omega), X_s(\omega))) \neq 0\},$$

which is an $(\mathcal{H}_s^t)$-stopping time by (A5.1) in [15]. In particular,

$$\{\omega \in \Omega \mid \mathbf{1}_{[\![0, \zeta[\![}(s, \omega) \mathbf{1}_{\Xi_t}(\omega)(A_s(\omega) - A_{s-}(\omega) - \varphi(X_{s-}(\omega), X_s(\omega))) = 0$$
$$\text{for } s \in [0, t[\}$$

$$= \{\omega \in \Omega \mid t < D_\Gamma(\omega)\} \in \mathcal{H}_t^t = \mathcal{F}_t.$$

Hence,

$$\{\omega \in \Xi_t \mid A_s(\omega) - A_{s-}(\omega) - \varphi(X_{s-}(\omega), X_s(\omega)) = 0 \text{ for } s \in {]}0, t{[} \cap {]}0, \zeta(\omega){[}\}$$
$$= \{\omega \in \Xi_t \mid A_s(\omega) - A_{s-}(\omega) - \varphi(X_{s-}(\omega), X_s(\omega)) = 0$$
$$\text{for } s \in [0, t{[} \cap [0, \zeta(\omega){[}\}$$
$$= \{\omega \in \Xi_t \mid \mathbf{1}_{[\![0, \zeta[\![}(s, \omega)(A_s(\omega) - A_{s-}(\omega) - \varphi(X_{s-}(\omega), X_s(\omega))) = 0$$
$$\text{for } s \in [0, t{[}\}$$
$$\in \mathcal{F}_t.$$

Therefore, $\widetilde{\Xi}_t \in \mathcal{F}_t$ and $\widetilde{\Xi} \in \mathcal{F}_\infty$. The proof for PrAF's on $[\![0, \zeta[\![$ is similar, so we omit it. □

The following theorem is a key to our extension of Nakao's operator $\Gamma$. Its proof is complicated by measurability issues, but the idea behind it is fairly transparent. We will use the convention $X_{0-}(\omega) := X_0(\omega)$.

THEOREM 2.18 (Dual PrAF). *Let $A$ be a finite rcll PrAF on $[\![0,\zeta[\![$ with defining sets $\Xi$, $\Xi_t$ admitting $m$-null set. Suppose that there is a Borel function $\varphi$ on $E \times E$ with $\varphi(x,x)=0$ for $x \in E$ such that $\varphi(X_{s-}(\omega), X_s(\omega)) = A_s(\omega) - A_{s-}(\omega)$ for all $s \in ]0,t[\cap]0,\zeta[$ and all $\omega \in \Xi_t$. Set*

(2.5)
$$\widehat{A}_t(\omega) := A_t(r_t(\omega)) + \varphi(X_t(\omega), X_{t-}(\omega))$$
*for $t \in [0, \zeta(\omega)[$ and $\widehat{A}_t(\omega) := 0$ for $t \in [\zeta(\omega), \infty[$.*

*$\widehat{A}$ is then an rcll PrAF on $[\![0,\zeta[\![$ admitting $m$-null set such that*

$$\widehat{A}_t = A_{t-} \circ r_t + \varphi(X_t, X_{t-}) \quad \text{and} \quad \widehat{A}_t - \widehat{A}_{t-} = \varphi(X_t, X_{t-})$$

*for all $t \in ]0,\zeta[$, $\mathbf{P}_m$-a.e.*

PROOF. Let $\Xi \in \mathcal{F}_\infty$, $\Xi_t \in \mathcal{F}_t^m$, $t > 0$ be the defining sets of $A$ admitting $m$-null set. We easily see $r_t^{-1}(\Xi_t) \cap \{t < \zeta\} \subset r_s^{-1}(\Xi_s) \cap \{s < \zeta\}$ for $s \in ]0,t[$ by use of Lemma 2.15(i) and $\theta_{t-s}\Xi_t \subset \Xi_s$.

Set $\widehat{\Xi}_t := r_t^{-1}(\Xi_t)$ for $t > 0$ and $\widehat{\Xi} := \bigcap_{t>0} \widehat{\Xi}_t$. We then see that $\widehat{\Xi} = \bigcap_{t>0, t\in\mathbb{Q}} \widehat{\Xi}_t$ by use of $r_t^{-1}(\Xi_t) \cap \{t \geq \zeta\} = \{t \geq \zeta\}$ and the monotonicity of $r_t^{-1}(\Xi_t) \cap \{t < \zeta\}$. Indeed, we have $\widehat{\Xi} \subset \bigcap_{t>0, t\in\mathbb{Q}} \widehat{\Xi}_t \subset (\widehat{\Xi}_s \cap \{s < \zeta\}) \cup \{t \geq \zeta\}$ for any $0 < s < t$ with $t \in \mathbb{Q}$. Taking the intersection over $t \in ]s, \infty[ \cap \mathbb{Q}$, we have $\widehat{\Xi} \subset \bigcap_{t>0, t\in\mathbb{Q}} \Xi_t \subset \widehat{\Xi}_s$ for all $s > 0$, which yields the assertion.

We prove $\theta_t \widehat{\Xi} \subset \widehat{\Xi}$ for each $t \geq 0$, in particular, $\theta_t \widehat{\Xi} \subset \theta_s \widehat{\Xi}$ and, equivalently, $\theta_s^{-1} \widehat{\Xi} \subset \theta_t^{-1} \widehat{\Xi}$ if $s \in [0,t]$. Suppose that $\omega \in \widehat{\Xi}$. Then, $r_{t+s}\omega \in \Xi_{t+s}$. If $t+s < \zeta(\omega)$, then $r_{t+s}\omega \in \Xi_s$, otherwise $r_{t+s}\omega = \omega_\Delta \in \Xi_s$. Hence, we have $r_s \theta_t \omega \in \Xi_s$, by Lemma 2.15(ii). Therefore, $r_s \theta_t \omega \in \Xi_s$ for all $s > 0$, which implies that $\theta_t \omega \in \widehat{\Xi}$.

Next, we prove $\theta_{t-s}(\widehat{\Xi}_t) \subset \widehat{\Xi}_s$ for $s \in ]0,t[$. Take $\omega \in \widehat{\Xi}_t$. Then, $r_s \theta_{t-s}\omega$ is pre-$s$-equivalent to $r_t \omega \in \Xi_t \subset \Xi_s$, by Lemma 2.15(ii) and hence $r_s \theta_{t-s}\omega \in \Xi_s$. Therefore, $\theta_{t-s}\omega \in \widehat{\Xi}_s$ for all $s \in ]0,t[$.

From $\Xi_t \subset \mathcal{F}_t^m$, we obtain $\widehat{\Xi}_t \in \mathcal{F}_t^m$, by Lemma 2.10. Since $(\widehat{\Xi}_t)^c = r_t^{-1}((\Xi_t)^c) = r_t^{-1}((\Xi_t)^c) \cap \{t < \zeta\}$ holds by noting $\omega_\Delta \in \Xi_t$, we have $\mathbf{P}_m((\widehat{\Xi}_t)^c) = \mathbf{P}_m((\widehat{\Xi})^c) = 0$.

By (2.2), $v \mapsto r_s(\omega)(v)$ is continuous at $v = s$. Hence, on $\widehat{\Xi}_t \cap \{t < \zeta\}$, we have $\varphi(X_{s-}, X_s) \circ r_s = \varphi(X_s, X_{s-}) \circ r_s = 0$, in particular, $A_s \circ r_s = A_{s-} \circ r_s$ for $s \in ]0,t[$.

The remainder of the proof is devoted to showing that $\widehat{A}$ is an rcll PrAF on $[0,\zeta[$ with defining sets $\widehat{\Xi}$, $\widehat{\Xi}_t$ such that on $\widehat{\Xi}_t \cap \{t < \zeta\}$, $\widehat{A}_s = A_{s-} \circ r_s + \varphi(X_s, X_{s-})$, $s \in ]0,t[$. First, note that for $\omega \in \widehat{\Xi}$, $|\widehat{A}_t(\omega)| < \infty$ for any $t \in ]0, \zeta(\omega)[$ because, by taking $T \in ]t, \zeta(\omega)[$, $r_T\omega \in \Xi_T$ implies that $r_t\omega \in \Xi_t$, hence $|A_{t-}(r_t\omega)| < \infty$. Moreover, for $\omega \in \widehat{\Xi}_t \cap \{t < \zeta\}$, we see that $r_s\omega \in \Xi_s \cap \{s < \zeta\}$ and $|A_{s-}(r_s\omega)| < \infty$ for all $0 < s < t$.

For two pre-$t$-equivalent paths $\omega, \omega' \in \Omega \cap \{t < \zeta\}$ with $t > 0$, we show that $\omega \in \widehat{\Xi}_t$ implies $\omega' \in \widehat{\Xi}_t$ and $\widehat{A}_s(\omega) = \widehat{A}_s(\omega')$ for $s \in [0, t[$. Recall that $\omega \in \widehat{\Xi}_t \cap \{t < \zeta\} \subset \widehat{\Xi}_s \cap \{s < \zeta\}$ for $s \in [0, t]$ and note that $\omega$ and $\omega'$ are $s$-equivalent for any $s \in [0, t[$. On the other hand, $s < \zeta(\omega)$ is equivalent to $s < \zeta(\omega')$ for any $s \in [0, t[$. We then see that $r_s \omega \in \Xi_s$ is $s$-equivalent to $r_s \omega'$ for any $s \in [0, t]$, which implies that $r_s \omega' \in \Xi_s$ for any $]0, t]$ and $A_{s-}(r_s \omega) = A_{s-}(r_s \omega')$ for any $s \in [0, t]$.

Fix $t > 0$. On $\widehat{\Xi}_t \cap \{t < \zeta\}$ and for any $p, q > 0$ with $p + q < t$, by Lemma 2.15,

$$\begin{aligned}
\widehat{A}_{p+q} &= A_{(p+q)-} \circ r_{p+q} + \varphi(X_{p+q}, X_{(p+q)-}) \\
&= (A_p + A_{q-} \circ \theta_p) \circ r_{p+q} + \varphi(X_{p+q}, X_{(p+q)-}) \\
&= A_p \circ r_{p+q} + A_{q-} \circ \theta_p \circ r_{p+q} + \varphi(X_{p+q}, X_{(p+q)-}) \\
&= (A_{p-} \circ r_{p+q} + \varphi(X_{p-}, X_p) \circ r_{p+q}) + A_{q-} \circ r_q + \varphi(X_{p+q}, X_{(p+q)-}) \\
&= (A_{p-} \circ r_p \circ \theta_q + \varphi(X_q, X_{q-})) + (\widehat{A}_q - \varphi(X_q, X_{q-})) \\
&\quad + \varphi(X_{p+q}, X_{(p+q)-}) \\
&= (\widehat{A}_p - \varphi(X_p, X_{p-})) \circ \theta_q + \widehat{A}_q + \varphi(X_{p+q}, X_{(p+q)-}) \\
&= \widehat{A}_p \circ \theta_q + \widehat{A}_q.
\end{aligned}$$

On $\widehat{\Xi}_t \cap \{t < \zeta\}$, again by Lemma 2.15 and (2.2), for any $s > 0$ and $u \in ]0, s[$,

$$\begin{aligned}
\widehat{A}_s - \widehat{A}_{s-u} &= \widehat{A}_u \circ \theta_{s-u} \\
&= (A_{u-} \circ r_u + \varphi(X_u, X_{u-})) \circ \theta_{s-u} \\
&= A_{u-} \circ r_u \circ \theta_{s-u} + \varphi(X_s, X_{s-}) \\
&= A_{u-} \circ r_s + \varphi(X_s, X_{s-}).
\end{aligned}$$

So,

$$\lim_{u \downarrow 0} (\widehat{A}_s - \widehat{A}_{s-u}) = \varphi(X_s, X_{s-}).$$

This shows that $\widehat{A}$ has left limit at $s \in ]0, t[$ and $\widehat{A}_s - \widehat{A}_{s-} = \varphi(X_s, X_{s-})$.

To show the right continuity of $\widehat{A}$ on $\widehat{\Xi}_t \cap \{t < \zeta\}$ at any $s \in ]0, t[$, note that for any $u \in ]0, t - s[$, by Lemma 2.15 and (2.2),

$$\begin{aligned}
\widehat{A}_{s+u} - \widehat{A}_s &= \widehat{A}_u \circ \theta_s \\
&= (A_{u-} \circ r_u + \varphi(X_u, X_{u-})) \circ \theta_s \\
&= A_{u-} \circ r_u \circ \theta_s + \varphi(X_{s+u}, X_{(s+u)-}) \\
&= A_{u-} \circ r_{s+u} + \varphi(X_{s+u}, X_{(s+u)-}).
\end{aligned}$$

Since $(A_v - A_{v-}) \circ r_{s+v} = \varphi(X_{v-}, X_v) \circ r_{s+v} = \varphi(X_s, X_{s-})$, while, by Lemma 2.15 and (2.2),

$$\begin{aligned}(A_v - A_{v-}) \circ r_{s+v} &= \lim_{u \downarrow 0}(A_v - A_{v-u}) \circ r_{s+v} \\ &= \lim_{u \downarrow 0} A_u \circ \theta_{v-u} \circ r_{s+v} \\ &= \lim_{u \downarrow 0} A_{u-} \circ r_{v+u} + \varphi(X_s, X_{s-}),\end{aligned}$$

we conclude that

$$\lim_{u \downarrow 0} A_{u-} \circ r_{s+u} = 0.$$

On the other hand, for any $s \geq 0$,

$$\begin{aligned}\lim_{u \downarrow 0} \varphi(X_{s+u}, X_{(s+u)-}) &= \lim_{u \downarrow 0} \varphi(X_{(v-u)-}, X_{v-u}) \circ r_{s+v} \\ &= \lim_{u \downarrow 0}(A_{v-u} - A_{(v-u)-}) \circ r_{s+v} \\ &= (A_{v-} - A_{v-}) \circ r_{s+v} = 0.\end{aligned}$$

Hence, we have, for $s > 0$,

$$\lim_{u \downarrow 0}(\widehat{A}_{s+u} - \widehat{A}_s) = 0.$$

In other words, $\widehat{A}$ is right-continuous at any $s \in ]0, t[$ on $\widehat{\Xi}_t \cap \{t < \zeta\}$. We also see that

$$\lim_{u < s, s \downarrow 0, u \downarrow 0}(\widehat{A}_{s+u} - \widehat{A}_s) = 0.$$

We can thus define the limit $\widehat{A}_0(\omega) := \lim_{s \downarrow 0} \widehat{A}_s(\omega)$ for $\omega \in \widehat{\Xi}_t \cap \{t < \zeta\}$ for any $t > 0$. We also see that $\widehat{A}_0(\omega) = \lim_{s \downarrow 0} \widehat{A}_{s-}(\omega)$ for $\omega \in \widehat{\Xi}_t \cap \{t < \zeta\}$ for any $t > 0$ because $\lim_{s \downarrow 0} \varphi(X_s, X_{s-}) = 0$. Next, we prove that $\widehat{A}_0(\omega) = 0$ for $\omega \in \widehat{\Xi}_t \cap \{t < \zeta\}$ for any $t > 0$. Take $\omega \in \widehat{\Xi}_t \cap \{t < \zeta\}$ for some fixed $t > 0$. It suffices to show that $\lim_{u \downarrow 0} \widehat{A}_{s-u}(\theta_u \omega) = \widehat{A}_s(\omega)$ for $s \in [0, t[$. Owing to Lemma 2.15(ii), we have

$$\begin{aligned}\widehat{A}_{s-u}(\theta_u \omega) &= A_{(s-u)-}(r_{s-u}\theta_u \omega) + \varphi(X_s(\omega), X_{s-}(\omega)) \\ &= A_{(s-u)-}(r_s \omega) + \varphi(X_s(\omega), X_{s-}(\omega)) \\ &= A_{s-u}(r_s \omega) - \varphi(X_u(\omega), X_{u-}(\omega)) + \varphi(X_s(\omega), X_{s-}(\omega)) \\ &= A_{s-u}(r_s \omega) - \widehat{A}_u(\omega) + \widehat{A}_{u-}(\omega) + \varphi(X_s(\omega), X_{s-}(\omega)) \\ &\to A_{s-}(r_s \omega) + \varphi(X_s(\omega), X_{s-}(\omega)) \qquad \text{as } u \downarrow 0 \\ &= \widehat{A}_s(\omega).\end{aligned}$$

The $\mathcal{F}_t^m$-measurability of $\widehat{A}_t$ is clear from (2.5). This proves the theorem. □

**3. Stochastic integral for Dirichlet processes.** The following fact will be used repeatedly in this section. Since a Hunt process is quasi-left continuous, for each fixed $t > 0$, we have $X_{t-} = X_t$, $\mathbf{P}_x$-a.s. for every $x \in E$.

Before embarking on the definition of our stochastic integral, we prepare the following lemma for later use.

LEMMA 3.1. *The following assertions hold.*

(i) *Let $\{G_n\}$ be an increasing sequence of finely open Borel sets. The following are then equivalent:*

  (a) *$\{G_n\}$ is a nest, that is, $\mathbf{P}_x(\lim_{n\to\infty} \sigma_{E\setminus G_n} \wedge \zeta = \zeta) = 1$ for q.e. $x \in E$;*
  (b) *$E = \bigcup_{n=1}^\infty G_n$ q.e.;*
  (c) *$\mathbf{P}_x(\lim_{n\to\infty} \sigma_{E\setminus G_n} = \infty) = 1$ for m-a.e. $x \in E$;*
  (d) *$\mathbf{P}_x(\lim_{n\to\infty} \sigma_{E\setminus G_n} = \infty) = 1$ for q.e. $x \in E$.*

*In particular, for an increasing sequence $\{F_n\}$ of closed sets, $\{F_n\}$ is an $\mathcal{E}$-nest if and only if $\mathbf{P}_x(\lim_{n\to\infty} \sigma_{E\setminus F_n} = \infty) = 1$ for m-a.e. $x \in E$.*

(ii) *For a function $f$ on $E$, $f \in \mathcal{F}_{\text{loc}}$ if and only if there exist an $\mathcal{E}$-nest $\{F_k\}$ of closed sets and $\{f_k \mid k \in \mathbb{N}\} \subset \mathcal{F}_b$ such that $f = f_k$ q.e. on $F_k$.*

PROOF. (i) For the implications (i)(a)$\Longleftrightarrow$(i)(b), see Theorem 4.6 in [11]. The implication (i)(d)$\Longrightarrow$(i)(a) is clear. Next, we show (i)(b)$\Longrightarrow$(i)(c). Since each $G_n$ is finely open, it is quasi-open by Theorem 4.6.1(i) in [7]. So, there exists a common nest $\{A_\ell\}$ of closed sets such that $(E \setminus G_n) \cap A_\ell$ is closed for all $n, \ell \in \mathbb{N}$. Set $\sigma := \lim_{n\to\infty} \sigma_{E\setminus G_n}$. We then have that for all $n \in \mathbb{N}$, $X_{\sigma_{E\setminus G_n}} \in E \setminus G_n$ $\mathbf{P}_x$-a.s. on $\{\sigma < \infty\}$ for q.e. $x \in E$. We have $\mathbf{P}_x(\lim_{\ell\to\infty} \sigma_{E\setminus A_\ell} = \infty) = 1$ q.e. $x \in E$. Since $\sigma(\omega) < \infty$ and $\lim_{\ell\to\infty} \sigma_{E\setminus A_\ell}(\omega) = \infty$ together imply $\sigma(\omega) < \sigma_{E\setminus A_{\ell_0}}(\omega)$ for some $\ell_0 = \ell_0(\omega) \in \mathbb{N}$, we have that there exists $\ell_0 \in \mathbb{N}$ such that $\sigma_{E\setminus G_n} < \sigma_{E\setminus A_\ell}$ for all $n > \ell \geq \ell_0$, $\mathbf{P}_x$-a.s. on $\{\sigma < \infty\}$ for q.e. $x \in E$. This means that

$$\mathbf{P}_x(\sigma < \infty) \leq \mathbf{P}_x\left(\varinjlim_{\ell\to\infty} \{X_{\sigma_{E\setminus G_n}} \in (E \setminus G_n) \cap A_\ell \text{ for all } n > \ell, \sigma < \infty\}\right)$$

$$\leq \varinjlim_{\ell\to\infty} \mathbf{P}_x(X_{\sigma_{E\setminus G_n}} \in (E \setminus G_\ell) \cap A_\ell \text{ for all } n > \ell, \sigma < \infty)$$

$$\leq \varinjlim_{\ell\to\infty} \mathbf{P}_x(X_\sigma \in (E \setminus G_\ell) \cap A_\ell, \sigma < \infty)$$

$$\leq \varinjlim_{\ell\to\infty} \mathbf{P}_x(X_\sigma \in E \setminus G_\ell, \sigma < \infty)$$

$$= \mathbf{P}_x\left(X_\sigma \in E \setminus \bigcup_{\ell=1}^\infty G_\ell, \sigma < \infty\right) = 0$$

for $m$-a.e. $x \in E$ because of the $\mathcal{E}$-polarity of $E \setminus \bigcup_{\ell=1}^\infty G_\ell$, where we use the quasi-left continuity of $X$ up to $\infty$ and the closedness of $(E \setminus G_\ell) \cap A_\ell$. The implication (i)(c) $\iff$ (i)(d) follows from the fact that $x \mapsto \mathbf{P}_x(\sigma < \infty)$ is the limit of a decreasing sequence of excessive functions and Lemma 4.1.7 in [7].

(i) The "if" part is clear by (i) because $\tau_{F_k} = \tau_{G_k}$, where $G_k$ is the fine interior of $F_k$. We only prove the "only if" part. Take $f \in \mathcal{F}_{\mathrm{loc}}$. There then exist $\{f_k \mid k \in \mathbb{N}\} \subset \mathcal{F}$ and an increasing sequence $\{G_k\}$ of finely open sets with $E = \bigcup_{k=1}^\infty G_k$ q.e. such that $f = f_k$ $m$-a.e. on $G_k$. We may take $f_k \in \mathcal{F}_b$ for each $k \in \mathbb{N}$ by replacing $f_k$ with $(-k) \vee f_k \wedge k$ and $G_k$ with $G_k \cap \{|f| < k\}$. Note that $f$ and $f_k$ are quasi-continuous, so $f = f_k$ q.e. on $G_k$. Taking an $\mathcal{E}$-quasi-closure $\overline{G_k}^{\mathcal{E}}$ of $G_k$, we have $f = f_k$ q.e. on $\overline{G_k}^{\mathcal{E}}$ (see [10] for the definition of $\mathcal{E}$-quasi-closure). Let $\{A_n\}$ be a common $\mathcal{E}$-nest of closed sets such that for each $k, n \in \mathbb{N}$, $\overline{G_k}^{\mathcal{E}} \cap A_n$ is closed. Set $F_k := \overline{G_k}^{\mathcal{E}} \cap A_k$. By (i), $\{G_k\}$ is a nest, hence $\overline{G_k}^{\mathcal{E}}$ is a nest of q.e. finely closed sets because of $\tau_{G_k} \leq \tau_{\overline{G_k}^{\mathcal{E}}}$. Here, we recognize $\overline{G_k}^{\mathcal{E}}$ as a finely closed Borel sets by deleting an $\mathcal{E}$-polar set. Since $\{A_n\}$ is a nest of closed sets, $\{F_k\}$ is also, that is, $\mathbf{P}_m(\lim_{k \to \infty} \tau_{F_k} \neq \zeta) = 0$. Therefore, $\{F_k\}$ is an $\mathcal{E}$-nest of closed sets. We easily see that for each $k \in \mathbb{N}$, $f = f_k$ q.e. on $F_k$. $\square$

Recall that any locally square-integrable MAF $M$ on $[\![0, \zeta[\![$ admits a jump function $\varphi$ on $E \times E$ with $\varphi(x,x) = 0$ for $x \in E$ such that $\Delta M_t = \varphi(X_{t-}, X_t)$ for $t \in \,]0, \zeta[$, $\mathbf{P}_m$-a.e. When $M \in \overset{\circ}{\mathcal{M}}$, we can strengthen this statement by replacing $]0, \zeta[$ with $]0, \infty[$ in view of Fukushima's decomposition and the combination of Theorem 5.2.1 and Lemma 5.6.3 in [7].

LEMMA 3.2. *Let $\phi$ be a Borel function on $E_\Delta \times E_\Delta$ satisfying $\phi(x,x) = 0$ for all $x \in E_\Delta$.*

(i) *Suppose that*

$$N(\mathbf{1}_{E \times E}(|\phi|^2 \wedge |\phi|))\mu_H \in S.$$

*There then exists a unique, purely discontinuous local MAF $K$ on $[\![0, \zeta[\![$ [i.e., $K \in (\mathcal{M}_{\mathrm{loc}}^1)^{[\![0, \zeta[\![}$] such that $K_t - K_{t-} = \phi(X_{t-}, X_t)$ for all $t \in [0, \zeta[$, $\mathbf{P}_x$-a.s. for q.e. $x \in E$.*

(ii) *If*

$$N(\mathbf{1}_{E \times E_\Delta}(|\phi|^2 \wedge |\phi|))\mu_H \in S,$$

*then $K$ can be taken to be a local MAF (i.e., $K \in \mathcal{M}_{\mathrm{loc}}^1$) and $K_t - K_{t-} = \phi(X_{t-}, X_t)$ for all $t \in [0, \infty[$, $\mathbf{P}_x$-a.s. for q.e. $x \in E$.*

PROOF. The proof of (ii) is similar to that of (i), so we only prove (i). By martingale theory (see, e.g. [9]), the hypothesis implies that the compensated

process

$$K_t^{(2)} := \sum_{0<s\leq t} \phi(X_{s-},X_s)\mathbf{1}_{\{|\phi(X_{s-},X_s)|>1\}}\mathbf{1}_{\{s<\zeta\}}$$
$$- \int_0^t \int_E N(X_s,dy)\phi(X_s,y)\mathbf{1}_{\{|\phi(X_s,y)|>1\}}\,dH_s$$

is a local MAF of $X$ of finite variation on $[\![0,\zeta[\![$ and

$$K_t^{(1)} := \lim_{\varepsilon \to 0}\Bigg(\sum_{0<s\leq t} \phi(X_{s-},X_s)\mathbf{1}_{\{\varepsilon<|\phi(X_{s-},X_s)|\leq 1\}}\mathbf{1}_{\{s<\zeta\}}$$
$$- \int_0^t \int_E N(X_s,dy)\phi(X_s,y)\mathbf{1}_{\{\varepsilon<|\phi(X_s,y)|\leq 1\}}\,dH_s\Bigg)$$

is a purely discontinuous locally square-integrable MAF of $X$ on $[\![0,\zeta[\![$. Thus, $K := K^{(1)} + K^{(2)}$ is a purely discontinuous MAF on $[\![0,\zeta[\![$ with jump function $\phi$. The uniqueness is clear from martingale theory. □

DEFINITION 3.3. Let $M$ be a local MAF on $[\![0,\zeta[\![$ [i.e., $M \in (\mathcal{M}_{\mathrm{loc}}^1)^{[\![0,\zeta[\![}]$ with jump function $\varphi$. Assume that for q.e. $x \in E$, $\mathbf{P}_x$-a.s.

$$\int_0^t \int_E (\widehat{\varphi}^2 \mathbf{1}_{\{|\widehat{\varphi}|\leq 1\}} + |\widehat{\varphi}|\mathbf{1}_{\{|\widehat{\varphi}|>1\}})(X_s,y)N(X_s,dy)\,dH_s < \infty \quad (3.1)$$
$$\text{for every } t<\zeta,$$

where $\widehat{\varphi}(x,y) := \varphi(x,y) + \varphi(y,x)$. Define, $\mathbf{P}_m$-a.e. on $[\![0,\zeta[\![$,

(3.2) $\Lambda(M)_t := -\frac{1}{2}(M_t + M_t \circ r_t + \varphi(X_t, X_{t-}) + K_t)$   for $t \in [0,\zeta[$,

where $K_t$ is the purely discontinuous local MAF on $[\![0,\zeta[\![$ with

(3.3) $\qquad K_t - K_{t-} = -\widehat{\varphi}(X_{t-},X_t)$   for every $t<\zeta$, $\mathbf{P}_x$-a.s.

for q.e. $x \in E$.

REMARK 3.4. (i) The condition (3.1) is nothing but $N(\mathbf{1}_{E\times E}(|\widehat{\varphi}|^2 \wedge |\widehat{\varphi}|))\mu_H \in S$. In particular, condition (3.1) is satisfied by the jump function of any element of $(\mathcal{M}_{\mathrm{loc}})^{[\![0,\zeta[\![}$.

(ii) It follows from Remark 2.13(iv) and Theorem 2.18 that $\Lambda(M)$ is a continuous PrAF admitting $m$-null set on $[\![0,\zeta[\![$. (This is because Remark 2.13(iv) and Theorem 2.18 imply that the process defined on $[0,\zeta[$ by $B_t := M_t \circ r_t + \varphi(X_t, X_{t-})$ is an rcll PrAF, with left-limit process $B_{t-} = B_t - \varphi(X_t, X_{t-})$. It follows that $\Lambda(M)$ is rcll on $[0,\zeta[$ and that $\Lambda(M)_{t-} = \Lambda(M)_t$ for all $t \in ]0,\zeta[$, $\mathbf{P}_m$-a.e.) Note that $-K_t := \sum_{s\leq t} \widehat{\varphi}(X_{s-},X_s)\mathbf{1}_{\{s<\zeta\}} -$ $\int_0^t \int_E \widehat{\varphi}(X_s, y) N(X_s, dy) \, dH_s$, $t < \zeta$, satisfies $K_t = K_t \circ r_t$ $\mathbf{P}_m$-a.e. on $\{t < \zeta\}$ for fixed $t > 0$. In view of Theorem 2.18, it is then clear from the definition that $\Lambda$ is a linear operator that maps local MAF's on $[\![0, \zeta[\![$ satisfying condition (3.1) into even CAF's on $[\![0, \zeta[\![$ admitting $m$-null set.

(iii) If $\{M^n, n \geq 1\}$ is a sequence of MAF's having finite energy and converging in probability to $M$, then it is easy to see that $M_t^n \circ r_t$, $\varphi^n(X_{t-}, X_t) = M_t^n - M_{t-}^n$ and $\varphi^n(X_t, X_{t-})$ converge to $M_t \circ r_t$, $\varphi(X_{t-}, X_t) = M_t - M_{t-}$ and $\varphi(X_t, X_{t-})$ in probability, respectively, under $\mathbf{P}_m$. Hence, we have that $\Lambda(M^n)_t$ converges to $\Lambda(M)_t$ in measure for each $t > 0$.

(iv) For $u \in \mathcal{F}$,
$$\Lambda(M^u)_t = -\tfrac{1}{2}(M_t^u + M_t^u \circ r_t + u(X_{t-}) - u(X_t)) = N_t^u,$$

$\mathbf{P}_m$-a.e. on $\{t < \zeta\}$, for each fixed $t \geq 0$. The first equality above is just the definition of $\Lambda(M^u)$, while the second follows by applying $r_t$ to both sides of (1.2) because $X_t = X_{t-}$ and $N^u \circ r_t = N_t^u$, $\mathbf{P}_m$-a.e. on $\{t < \zeta\}$. (The last property is proved in [4], Theorem 2.1, when $X$ is a diffusion, but the same proof works for general symmetric Markov process $X$.) Since both $\Lambda(M^u)_t$ and $N_t^u$ are continuous in $t$, we even have, $\mathbf{P}_m$-a.e.,
$$\Lambda(M^u)_t = N_t^u \qquad \text{for all } t < \zeta.$$

We are going to show that $\Lambda(M)$ defined above coincides on $[0, \zeta[$ with $\Gamma(M)$ defined in (1.5) by Nakao when $M$ is an MAF of finite energy. An AF $Z$ is called *even* (resp., *odd*) if and only if $Z_t \circ r_t = Z_t$ (resp., $Z_t \circ r_t = -Z_t$), $\mathbf{P}_m$-a.e. on $\{t < \zeta\}$ for each $t > 0$. For an rcll process $Z$ with $Z_0 = 0$ and $T > 0$, we define
$$R_T Z_t := (R_T Z)_t := Z_{T-} - Z_{(T-t)-} \qquad \text{for } 0 \leq t \leq T,$$

with the convention $Z_{0-} = Z_0 = 0$. Note that $R_T Z_t$ so defined is an rcll process in $t \in [0, T]$.

LEMMA 3.5. *Suppose that $Z$ is an rcll PrAF. Then, $\mathbf{P}_m$-a.e. on $\{T < \zeta\}$,*

(3.4) $\quad R_T Z_t = \begin{cases} Z_t \circ r_T, & \text{if } Z \text{ is even,} \\ -Z_t \circ r_T, & \text{if } Z \text{ is odd,} \end{cases} \qquad \text{for every } t \in [0, T].$

PROOF. Let $Z$ be an rcll PrAF and let $T > 0$. By Lemma 2.15,
$$Z_t \circ r_T = (Z_T - Z_{T-t} \circ \theta_t) \circ r_T = Z_T \circ r_T - Z_{T-t} \circ r_{T-t}$$
(3.5)
$$\text{for all } t < T.$$

When $Z$ is even,
$$Z_t \circ r_T = Z_T - Z_{T-t} = Z_{T-} - Z_{(T-t)-} = R_T Z_t,$$

$\mathbf{P}_m$-a.e. on $\{T < \zeta\}$ for each fixed $0 \leq t < T$. Since both sides are right-continuous in $t \in [0, T[$, we have, $\mathbf{P}_m$-a.e., $R_T Z_t = Z_t \circ r_T$ for every $t \in [0, T]$. When $Z$ is an odd AF of $Z$, (3.4) can be proven similarly. □

THEOREM 3.6. *For an MAF $M$ of finite energy, $\Lambda(M)$ defined above coincides on $[\![0, \zeta[\![$ with $\Gamma(M)$ defined in (1.5), $\mathbf{P}_m$-a.e.*

PROOF. For $u \in \mathcal{F}$ and $0 < t < T$, since $N^u$ is an even CAF, by Lemma 3.5,

$$\begin{aligned}(M_t^u + 2N_t^u) \circ r_T &= (u(X_t) - u(X_0) + N_t^u) \circ r_T \\ &= u(X_{(T-t)-}) - u(X_{T-}) + N_{T-}^u - N_{(T-t)-}^u \\ &= M_{(T-t)-}^u - M_{T-}^u \\ &= -R_T M_t^u.\end{aligned}$$

Since both $(M_t^u + 2N_t^u) \circ r_T$ and $R_T M_t^u$ are right-continuous in $t$, we have, $\mathbf{P}_m$-a.e. on $\{T < \zeta\}$,

(3.6) $\qquad R_T M_t^u = -(M_t^u + 2N_t^u) \circ r_T \qquad$ for every $t \in [0, T]$.

For $u \in \mathcal{D}(\mathcal{L}) \subset \mathcal{F}$ and $v \in \mathcal{F}_b$, define $M_t = \int_0^t v(X_{s-}) \, dM_s^u$, which is an MAF of finite energy. Note that, since $u \in \mathcal{D}(\mathcal{L})$, $N_t^u = \int_0^t \mathcal{L}u(X_s) \, ds$ is a continuous process of finite variation. For each fixed $0 < t < T$ and $n \geq 1$, define $t_i = it/n$ and $s_i = T - t + t_i$. Using the standard Riemann sum approximation of the Itô integral and of the covariance process $[M^v, M^u]$, we have, $\mathbf{P}_m$-a.e. on $\{T < \zeta\}$,

$$\begin{aligned}& M_T - M_{T-t} + [M^v, M^u]_T - [M^v, M^u]_{T-t} \\ &= \lim_{n \to \infty} \left( \sum_{i=0}^{n-1} v(X_{s_i})(M_{s_{i+1}}^u - M_{s_i}^u) + (M_{s_{i+1}}^v - M_{s_i}^v)(M_{s_{i+1}}^u - M_{s_i}^u) \right) \\ &= \lim_{n \to \infty} \left( \sum_{i=0}^{n-1} v(X_{s_{i+1}})(M_{s_{i+1}}^u - M_{s_i}^u) - (N_{s_{i+1}}^v - N_{s_i}^v)(M_{s_{i+1}}^u - M_{s_i}^u) \right) \\ &= \lim_{n \to \infty} \sum_{i=0}^{n-1} v(X_{s_{i+1}})(M_{s_{i+1}}^u - M_{s_i}^u) \\ &= \lim_{n \to \infty} \sum_{i=0}^{n-1} v(X_{T-t+t_i})(R_T M_{t-t_i}^u - R_T M_{t-t_{i+1}}^u) \\ &= \lim_{n \to \infty} \left( \sum_{i=0}^{n-1} v(X_{t-t_{i+1}})(M_{t-t_{i+1}}^u - M_{t-t_i}^u + 2N_{t-t_{i+1}}^u - 2N_{t-t_i}^u) \right) \circ r_T \\ &= -\left( \int_0^t v(X_{s-}) \, d(M_s^u + 2N_s^u) \right) \circ r_T,\end{aligned}$$

where in the third equality, we used the fact that $N^u$ has zero energy, while in the second to the last equality, we used (3.6). Note that the stochastic integral involving $N^u$ in the last equality is just the Lebesgue–Stieltjes integral since $N^u$ is of finite variation. Also, note that $X_t = X_{t-}$, $\mathbf{P}_m$-a.e., for each fixed $t > 0$. So, we have, for each fixed $t < T$, $\mathbf{P}_m$-a.e. on $\{T < \zeta\}$,

$$R_T M_t + R_T [M^v, M^u]_t = -\left(\int_0^t v(X_{s-}) \, d(M_s^u + 2N_s^u)\right) \circ r_T.$$

Since both sides are right-continuous in $t \in [0, T]$, we have, $\mathbf{P}_m$-a.e. on $\{T < \zeta\}$,

(3.7)
$$R_T M_t + R_T [M^v, M^u]_t = -\left(\int_0^t v(X_{s-}) d(M_s^u + 2N_s^u)\right) \circ r_T$$
for every $t \in [0, T]$.

By [14], Theorem 3.1 and (1.7),

$$\int_0^t v(X_{s-}) \, dN_s^u = \int_0^t v(X_s) \, dN_s^u = \Gamma(M)_t - \tfrac{1}{2}\langle M^{v,c} + M^{v,j}, M^{u,c} + M^{u,j}\rangle_t.$$

It follows that $\mathbf{P}_m$-a.e. on $\{T < \zeta\}$,

$$\begin{aligned} R_T M_t &+ R_T [M^v, M^u]_t \\ &= -(M_t + 2\Gamma(M)_t - \langle M^{v,c} + M^{v,j}, M^{u,c} + M^{u,j}\rangle_t) \circ r_T \\ &= -(M_t + 2\Gamma(M)_t - \langle M^{v,c}, M^{u,c}\rangle_t - \langle M^{v,j}, M^{u,j}\rangle_t) \circ r_T \end{aligned}$$
for all $t \leq T$.

Recall that

$$\begin{aligned}[M^v, M^u]_t &= \langle M^{v,c}, M^{u,c}\rangle_t + \sum_{s \leq t}(M_s^v - M_{s-}^v)(M_s^u - M_{s-}^u) \\ &= \langle M^{v,c}, M^{u,c}\rangle_t + \sum_{s \leq t}(v(X_s) - v(X_{s-}))(u(X_s) - u(X_{s-})).\end{aligned}$$

Taking $t = T$ and noting that both $\Gamma(M)$ and $\langle M^{v,c}, M^{u,c}\rangle$ are continuous even AF's, we have, from above, that, $\mathbf{P}_m$-a.e. on $\{t < \zeta\}$,

(3.8)   $\Gamma(M)_t = -\tfrac{1}{2}(M_t + M_t \circ r_t + v(X_t)(u(X_{t-}) - u(X_t)) + K_t),$

where

$$K_t = \sum_{s \leq t}(v(X_s) - v(X_{s-}))(u(X_s) - u(X_{s-})) - \langle M^{v,j}, M^{u,j}\rangle_t$$

is the purely discontinuous MAF with $K_t - K_{t-} = (v(X_t) - v(X_{t-}))(u(X_t) - u(X_{t-}))$. Note that the right-hand side of (3.8) is right-continuous on $[0, \zeta[$, $\mathbf{P}_m$-a.e. [cf. Remark 3.4(ii)]. Also, observe that $M_t - M_{t-} = \varphi(X_{t-}, X_t)$, where $\varphi(x,y) = v(x)(u(y) - u(x))$, and that

$$K_t - K_{t-} = -\varphi(X_{t-}, X_t) - \varphi(X_t, X_{t-}).$$

This shows that $\Gamma(M)_t = \Lambda(M)_t$, $\mathbf{P}_m$-a.e. on $\{t < \zeta\}$ for each fixed $t \geq 0$. Since both processes are continuous in $t \in [0, \zeta[$, we have, $\mathbf{P}_m$-a.e.,

$$\Gamma(M) = \Lambda(M) \qquad \text{on } [0, \zeta[$$

for an MAF $M$ of the form $M_t = \int_0^t v(X_{s-}) \, dM_s^u$ with $u \in \mathcal{D}(\mathcal{L})$ and $v \in \mathcal{F}_b$. By Lemma 5.4.5 in [6], such MAF's form a dense subset in the space of MAF's having finite energy. Thus, by Lemma 3.1 in Nakao [14] and Remark 3.4(iii), we have, for a general MAF $M$ of finite energy, $\Gamma(M)_t = \Lambda(M)_t$ $\mathbf{P}_m$-a.e. on $\{t < \zeta\}$ for every fixed $t \geq 0$. Since both processes are continuous in $t \in [0, \zeta[$, it follows that $\Gamma(M) = \Lambda(M)$ on $[\![0, \zeta[\![$, $\mathbf{P}_m$-a.e. □

THEOREM 3.7. *Let $M$ be a locally square-integrable MAF on $[\![0, \zeta[\![$ with jump function $\varphi$. Suppose that $\varphi$ satisfies condition* (3.1). *Then, for every $t > 0$,*

$$(3.9) \qquad \lim_{n \to \infty} \sum_{\ell=0}^{n-1} (\Lambda(M)_{(\ell+1)t/n} - \Lambda(M)_{\ell t/n})^2 = 0,$$

*where the convergence is in $\mathbf{P}_{gm}$-measure on $\{t < \zeta\}$ for any $g \in L^1(E;m)$ with $0 < g \leq 1$ $m$-a.e.*

PROOF. By (1.5) and Theorem 3.6, (3.9) clearly holds when $M$ is an MAF of finite energy. For a locally square-integrable MAF $M$ on $[\![0, \zeta[\![$, there is an $\mathcal{E}$-nest $\{F_k\}$ of closed sets such that $\mathbf{1}_{F_k} * M \in \overset{\circ}{\mathcal{M}}$ for each $k \geq 1$ in view of the proof of Proposition 2.8 and so (3.9) holds with $\mathbf{1}_{F_k} * M$ in place of $M$. For each fixed $k \geq 1$,

$$\Lambda(M)_t = \Lambda(\mathbf{1}_{F_k} * M)_t - \tfrac{1}{2} K_t^k, \qquad \mathbf{P}_m\text{-a.e. on } [0, \tau_{F_k}[,$$

where $K_t^k$ is a purely discontinuous local MAF on $[\![0, \zeta[\![$ with

$$K_t^k - K_{t-}^k = \mathbf{1}_{F_k^c}(X_{t-})\varphi(X_{t-}, X_t) + \mathbf{1}_{F_k^c}(X_t)\varphi(X_t, X_{t-}) \qquad \text{for } t < \zeta.$$

Since $\mathbf{1}_{F_k} * M \in \overset{\circ}{\mathcal{M}}$, we have

$$\int_E N(\mathbf{1}_{F_k \times E} \varphi^2) \, d\mu_H = \int_E N(\mathbf{1}_{E \times F_k} \varphi^2) \, d\mu_H < \infty.$$

Consequently, by Lemma 3.2, we have the existence of a purely discontinuous local MAF on $[\![0, \zeta[\![$ with jumps given by $\mathbf{1}_{F_k}(X_{t-})\varphi(X_{t-}, X_t) + \mathbf{1}_{F_k}(X_t)\varphi(X_t, X_{t-})$, $t < \zeta$. So, we obtain the existence of such $K_t^k$. Since the square bracket of $K^k$ is given by $\sum_{s\leq t}\mathbf{1}_{F_k^c}(X_{s-})\varphi^2(X_{s-},X_s) + \mathbf{1}_{F_k^c}(X_s)\varphi^2(X_s,X_{s-})$ and it vanishes at $t < \tau_{F_k}$, we have, for each fixed $t > 0$,

$$\lim_{n\to\infty}\sum_{\ell=0}^{n-1}(\Lambda(M)_{(\ell+1)t/n} - \Lambda(M)_{\ell t/n})^2 = 0 \qquad \text{in } \mathbf{P}_{gm}\text{-measure on } \{t < \tau_{F_k}\}.$$

Passing to the limit as $k\uparrow\infty$ establishes (3.9). $\square$

We are now in a position to define stochastic integrals against $\Lambda(M)$ as integrator. Note that, for $f\in\mathcal{F}_{\text{loc}}$, $M^{f,c}$ is well defined as a continuous MAF on $[\![0,\zeta[\![$ of locally finite energy (see Theorem 8.2 in [9]). Moreover, for $f\in\mathcal{F}_{\text{loc}}$ and a locally square-integrable MAF $M$ on $[\![0,\zeta[\![$,

$$t\mapsto (f*M)_t := \int_0^t f(X_{s-})\,dM_s$$

is a locally square-integrable MAF on $[\![0,\zeta[\![$.

DEFINITION 3.8 (*Stochastic integral*). Suppose that $M$ is a locally square-integrable MAF on $[\![0,\zeta[\![$ and $f\in\mathcal{F}_{\text{loc}}$. Let $\varphi: E\times E\to \mathbb{R}$ be a jump function for $M$ and assume that $\varphi$ satisfies condition (3.1). Define, $\mathbf{P}_m$-a.e. on $[\![0,\zeta[\![$,

$$
\begin{aligned}
\int_0^t f(X_{s-})\,d\Lambda(M)_s \\
:= \Lambda(f*M)_t - \tfrac{1}{2}\langle M^{f,c},M^c\rangle_t \\
+ \tfrac{1}{2}\int_0^t\int_E (f(y)-f(X_s))\varphi(y,X_s)N(X_s,dy)\,dH_s,
\end{aligned}
$$
(3.10)

whenever $\Lambda(f*M)$ is well defined and the third term in the right-hand side of (3.10) is absolutely convergent.

REMARK 3.9. (i) Under the above condition, the stochastic integral is clearly well defined on $[\![0,\zeta[\![$ under $\mathbf{P}_m$ and is a PrAF of $X$ admitting $m$-null set.

(ii) Here are some sufficient conditions for every term on the right-hand side of (3.10) to be well defined. In addition to the conditions in Definition 3.8, we assume that, $\mathbf{P}_m$-a.e.,

(3.11) $\quad \int_0^t\int_E (f(X_s)-f(y))^2 N(X_s,dy)\,dH_s < \infty \qquad$ for every $t < \zeta$

and that

(3.12) $\quad \int_0^t\int_E \varphi(y,X_s)^2 N(X_s,dy)\,dH_s < \infty \qquad$ for every $t < \zeta$.

The first and third terms on the right-hand side of (3.10) are then well defined. This is because $N(\mathbf{1}_{E\times E}|\widehat{\varphi}|)\mu_H \in S$ implies that $N(\mathbf{1}_{E\times E}|f\widehat{\varphi}|)\mu_H \in S$, and

$$f(x)\varphi(x,y) + f(y)\varphi(y,x) = f(x)\hat{\varphi}(x,y) + (f(y) - f(x))\varphi(y,x),$$

so $\Lambda(f*M)$ is well defined on $[0,\zeta[$ in view of the condition (3.1) for $f*M$, (3.11) and (3.12). Condition (3.11) is satisfied when $f$ is a bounded function in $\mathcal{F}_{\mathrm{loc}}$ or $f \in \mathcal{F}$. This is because, when $f \in \mathcal{F}$, the left-hand side of (3.11) is just $\langle M^{f,j}\rangle_t$. When $f$ is a bounded function in $\mathcal{F}_{\mathrm{loc}}$, by Lemma 3.1(i), there exist a nest $\{F_n \mid n \in \mathbb{N}\}$ of closed sets and a sequence of functions $\{f_n \mid n \in \mathbb{N}\} \subset \mathcal{F}_b$ such that $f = f_n$ q.e. on $F_n$ for every $n \geq 1$. Note that for each $n \geq 1$, $M^{f_n,d}$ is a square-integrable, purely discontinuous martingale and

$$M_t^{f_n,d} - M_{t-}^{f_n,d} = f_n(X_t) - f_n(X_{t-}).$$

So, $t \mapsto \sum_{s\leq t}(f_n(X_s) - f_n(X_{s-}))^2$ is $\mathbf{P}_x$-integrable for q.e. $x \in E$. Since $f$ is bounded, we have, for each $n \geq 1$, that

$$t \mapsto \sum_{s\leq t\wedge\tau_{F_n}} (f(X_s) - f(X_{s-}))^2$$

$$= \sum_{s<t\wedge\tau_{F_n}} (f(X_s) - f(X_{s-}))^2 + (f(X_{t\wedge\tau_{F_n}}) - f(X_{t\wedge\tau_{F_n}-}))^2$$

$$= \sum_{s<t\wedge\tau_{F_n}} (f_n(X_s) - f_n(X_{s-}))^2 + (f(X_{t\wedge\tau_{F_n}}) - f(X_{t\wedge\tau_{F_n}-}))^2$$

is an increasing process and is $\mathbf{P}_x$-integrable for each fixed $t \geq 0$ for q.e. $x \in E$. Similarly, $A_t := \sum_{s\leq t}(f(X_s) - f(X_{s-}))^2$ is *locally integrable* in the sense of Definition 5.18 in [9]. Indeed, for a stopping time $T_n := \inf\{t > 0 \mid A_t > n\}$, $A_{T_n} = A_{T_n-} + (f(X_{T_n}) - f(X_{T_n-}))^2$ is bounded, hence $\mathbf{P}_x$-integrable for q.e. $x \in E$. Note that the dual predictable projection of $A_t$ is nothing but $\int_0^t \int_{E_\Delta} (f(X_s) - f(y))^2 N(X_s, dy) dH_s$. The dual predictable projection of $\sum_{s\leq t\wedge\tau_{F_n}} (f(X_s) - f(X_{s-}))^2$ is then given by $\int_0^{t\wedge\tau_{F_n}} \int_{E_\Delta} (f(X_s) - f(y))^2 N(X_s, dy) dH_s$ from Corollary 5.24 in [9], which is $\mathbf{P}_x$-integrable for q.e. $x \in E$. This implies that (3.11) holds for every $t < \tau_{F_n}$. Therefore, (3.11) holds for every $t < \zeta$.

Condition (3.12) is satisfied when $M^d$ is $\mathbf{P}_m$-square-integrable. Indeed,

$$\mathbf{E}_m\left[\sum_{s\leq t}\varphi^2(X_s, X_{s-}) : t < \zeta\right] = \mathbf{E}_m[[M^d]_t \circ r_t : t < \zeta]$$

$$= \mathbf{E}_m[[M^d]_t : t < \zeta] < \infty.$$

Corollary 4.5 in [8] then tells us that

$$\lim_{t\to 0}\frac{1}{t}\mathbf{E}_m\left[\sum_{s\leq t}\varphi^2(X_s,X_{s-}):t<\zeta\right]=\lim_{t\to 0}\frac{1}{t}\mathbf{E}_m\left[\sum_{s\leq t}\varphi^2(X_s,X_{s-})\right],$$

which implies that

$$\mathbf{E}_m\left[\int_0^t\int_E\varphi(y,X_s)^2N(X_s,dy)\,dH_s\right]<\infty$$

for all $t>0$ by way of its subadditivity. Hence, we obtain (3.12).

(iii) Suppose that $f\in\mathcal{F}$. Let $K_t$ be a purely discontinuous local MAF on $[\![0,\zeta[\![$ with $K_t-K_{t-}=-\varphi(X_{t-},X_t)-\varphi(X_t,X_{t-})$ on $]0,\zeta[$. Then,

$$\langle M^{f,j},M^j+K\rangle_t=-\int_0^t\int_E(f(y)-f(X_s))\varphi(y,X_s)N(X_s,dy)\,dH_s.$$

In this case, (3.10) can be rewritten as

$$(3.13)\quad \int_0^t f(X_{s-})\,d\Lambda(M)_s=\Lambda(f*M)_t-\tfrac{1}{2}\langle M^{f,c}+M^{f,j},M^c+M^j+K\rangle_t$$

on $[0,\zeta[$. So, when $M=M^u$ for some $u\in\mathcal{F}$ and $f\in\mathcal{F}\cap L^2(E;\mu_{\langle u\rangle})$, $\int_0^t f(X_{s-})\,d\Lambda(M)_s$ on $[0,\zeta[$ is just the $\int_0^t f(X_s)\circ d\Gamma(M)_s$ defined by (1.7). This shows that the stochastic integral given in Definition 3.8 extends Nakao's definition (1.7) of stochastic integral first introduced in [14].

THEOREM 3.10. *The stochastic integral in* (3.10) *is well defined. That is, if $M$ and $\widetilde{M}$ are two locally square-integrable MAF's on $[\![0,\zeta[\![$ such that all conditions in Definition* 3.3 *for $M$ and $\widetilde{M}$ are satisfied and $\Lambda(M)\equiv\Lambda(\widetilde{M})$ on $[\![0,\zeta[\![$, then for every $f\in\mathcal{F}_{\mathrm{loc}}$ for which $\int_0^t f(X_{s-})d\Lambda(M)_s$ and $\int_0^t f(X_{s-})\,d\Lambda(\widetilde{M})_s$ are well defined, we have,* $\mathbf{P}_m$-*a.e.,*

$$\int_0^t f(X_{s-})\,d\Lambda(M)_s=\int_0^t f(X_{s-})\,d\Lambda(\widetilde{M})_s\qquad on\ [0,\zeta[.$$

PROOF. It is equivalent to show that

$$\int_0^t f(X_{s-})\,d\Lambda(M-\widetilde{M})_s=0\qquad \text{on }[0,\zeta[.$$

By taking $M$ to be $M-\widetilde{M}$, we may and will assume that $\widetilde{M}=0$. Moreover, a localization argument allows us to assume that $f$ is bounded. Let $\varphi:E\times E\to\mathbb{R}$ be a jump function for $M$. Let $K_t$ be the purely discontinuous local MAF on $[\![0,\zeta[\![$ with

$$K_t-K_{t-}=-\varphi(X_{t-},X_t)-\varphi(X_t,X_{t-})\qquad\text{for }t<\zeta.$$

Since $\Lambda(M) \equiv 0$, we have

(3.14) $\quad M_t + M_t \circ r_t + \varphi(X_t, X_{t-}) + K_t = 0 \qquad \text{on } [0, \zeta[.$

Thus, by (3.5) and (3.14), on $\{T < \zeta\}$,

$$\begin{aligned}(3.15) \quad M_t \circ r_T &= M_T \circ r_T - M_{T-t} \circ r_{T-t} \\ &= -M_T - K_T + M_{T-t} + K_{T-t} \\ &\quad - \varphi(X_T, X_{T-}) + \varphi(X_{T-t}, X_{(T-t)-})\end{aligned}$$

for every $t \in [0, T]$. Using the standard Riemann sum approximation and (3.15), we have, for $f \in \mathcal{F}$,

$$\begin{aligned}(f * M)_t \circ r_t + f(X_t)\varphi(X_t, X_{t-}) \\ = -(f * M)_t - (f * K)_t - [M^f, M + K]_t \\ = -(f * M)_t - (f * K)_t - \langle M^{f,c}, M^c \rangle_t \\ + \sum_{s \le t}(f(X_s) - f(X_{s-}))\varphi(X_s, X_{s-})\end{aligned}$$

$\mathbf{P}_m$-a.e. on $\{t < \zeta\}$ for each fixed $t \ge 0$. Consequently, we have, for $f \in \mathcal{F}_{\text{loc}}$, $\mathbf{P}_m$-a.e. for all $t \in [0, \zeta[$,

$$\begin{aligned}(3.16) \quad & (f * M)_t \circ r_t + f(X_t)\varphi(X_t, X_{t-}) \\ &= -(f * M)_t - (f * K)_t - \langle M^{f,c}, M^c \rangle_t \\ &\quad + \sum_{s \le t}(f(X_s) - f(X_{s-}))\varphi(X_s, X_{s-}),\end{aligned}$$

since both sides are right-continuous in $t \in [0, \zeta[$. Let $\widetilde{K}$ be the purely discontinuous local MAF on $[\![0, \zeta[\![$ with

$$\widetilde{K}_t - \widetilde{K}_{t-} = -f(X_{t-})\varphi(X_{t-}, X_t) - f(X_t)\varphi(X_t, X_{t-}) \qquad \text{for all } t \in [0, \zeta[.$$

Then, for $f \in \mathcal{F}_{\text{loc}}$, we have, by (3.16),

$$\begin{aligned}\Lambda(f * M)_t &= -\tfrac{1}{2}((f * M)_t + (f * M) \circ r_t + f(X_t)\varphi(X_t, X_{t-}) + \widetilde{K}_t) \\ &= \tfrac{1}{2}\left(\int_0^t f(X_{s-})\, dK_s + \langle M^{f,c}, M^c \rangle_t \right.\\ &\quad \left. - \sum_{s \le t}(f(X_s) - f(X_{s-}))\varphi(X_s, X_{s-}) - \widetilde{K}_t\right).\end{aligned}$$

Thus,

$$\int_0^t f(X_{s-})\, d\Lambda(M)_s$$

$$= \Lambda(f * M)_t - \tfrac{1}{2}\langle M^{f,c}, M^c\rangle_t$$
$$+ \tfrac{1}{2}\int_0^t \int_E (f(y) - f(X_s))\varphi(y, X_s) N(X_s, dy)\, dH_s$$
$$= \tfrac{1}{2}\int_0^t f(X_{s-})\, dK_s - \tfrac{1}{2}\sum_{s\leq t}(f(X_s) - f(X_{s-}))\varphi(X_s, X_{s-}) - \tfrac{1}{2}\widetilde{K}_t$$
$$+ \tfrac{1}{2}\int_0^t \int_E (f(y) - f(X_s))\varphi(y, X_s) N(X_s, dy)\, dH_s.$$

Note that

$$(3.17)\quad \widetilde{K}_t = -\sum_{s\leq t}(f(X_{s-})\varphi(X_{s-}, X_s) + f(X_s)\varphi(X_s, X_{s-}))$$
$$+ \int_0^t \int_E (f(X_s)\varphi(X_s, y) + f(y)\varphi(y, X_s)) N(X_s, dy)\, dH_s$$

and that

$$(3.18)\quad K_t = \lim_{\varepsilon\to 0}\left(-\sum_{s\leq t}(\widehat{\varphi}\mathbf{1}_{\{|\widehat{\varphi}|>\varepsilon\}})(X_{s-}, X_s) + (N(\widehat{\varphi}\mathbf{1}_{\{|\widehat{\varphi}|>\varepsilon\}}) * H)_t\right),$$

where $\widehat{\varphi}(x, y) := \varphi(x, y) + \varphi(y, x)$. It follows that

$$\int_0^t f(X_{s-})\, d\Lambda(M)_s = 0 \qquad \text{for all } t < \zeta,$$

$\mathbf{P}_m$-a.e. □

REMARK 3.11. The above proof actually shows that if $\Lambda(M) = \Lambda(\widetilde{M})$ on $[0, T] \cap [0, \zeta[$, then, $\mathbf{P}_m$-a.e.,

$$\int_0^t f(X_{s-})\, d\Lambda(M)_s = \int_0^t f(X_{s-})\, d\Lambda(\widetilde{M})_s \qquad \text{on } [0, T] \cap [0, \zeta[.$$

**4. Further study of the stochastic integral.**

THEOREM 4.1. *Suppose that $f \in \mathcal{F}_{\mathrm{loc}}$ and that $M$ is a locally square-integrable MAF on $[\![0, \zeta[\![$ satisfying (3.1) such that $\Lambda(M)$ is a continuous process $A$ of finite variation on $[\![0, \zeta[\![$. Assume that the stochastic integral $t \mapsto \int_0^t f(X_{s-})\, d\Lambda(M)_s$ is well defined. Then, $\mathbf{P}_m$-a.e.,*

$$\int_0^t f(X_{s-})\, d\Lambda(M)_s = \int_0^t f(X_s)\, dA_s \qquad \text{on } [0, \zeta[,$$

*where the integral on the right-hand side is the Lebesgue–Stieltjes integral.*

PROOF. Let $\varphi: E \times E \to \mathbb{R}$ be a Borel function with $\varphi(x,x) = 0$ for $x \in E$ such that $\varphi(X_{t-}, X_t) = M_t - M_{t-}$ for $t \in [0, \zeta[$, $\mathbf{P}_m$-a.e. Let $K_t$ be the purely discontinuous local MAF on $[\![0, \zeta[\![$ with

$$K_t - K_{t-} = -\varphi(X_{t-}, X_t) - \varphi(X_t, X_{t-}) \qquad \text{for } t \in ]0, \zeta[.$$

Since $\Lambda(M) = A$ on $[0, \zeta[$,

$$M_t \circ r_t + \varphi(X_t, X_{t-}) = -M_t - K_t - 2A_t \qquad \text{for all } t \in [0, \zeta[.$$

Thus, by (3.5), for every $T > t > 0$, on $\{T < \zeta\}$,

$$\begin{aligned}(4.1)\quad M_t \circ r_T &= -M_T - K_T - 2A_T + M_{T-t} + K_{T-t} \\ &\quad + 2A_{T-t} - \varphi(X_T, X_{T-}) + \varphi(X_{T-t}, X_{(T-t)-}).\end{aligned}$$

Now, fix $f \in \mathcal{F}_{\mathrm{loc}}$; as before, we may assume, without loss of generality, that $f$ is bounded. Using the standard Riemann sum approximation, we obtain, on $\{t < \zeta\}$,

$$\begin{aligned}(f * M)_t &\circ r_t + f(X_t)\varphi(X_t, X_{t-}) \\ &= -(f*M)_t - (f*K)_t - 2(f*A)_t - [M^f, M+K+2A]_t \\ &= -(f*M)_t - (f*K)_t - 2(f*A)_t - \langle M^{f,c}, M^c \rangle_t \\ &\quad + \sum_{s \leq t}(f(X_s) - f(X_{s-}))\varphi(X_s, X_{s-}).\end{aligned}$$

Consequently, we have, $\mathbf{P}_m$-a.e. for all $t \in [0, \zeta[$,

$$\begin{aligned}(4.2)\quad (f*M)_t &\circ r_t + f(X_t)\varphi(X_t, X_{t-}) \\ &= -(f*M)_t - (f*K)_t - 2(f*A)_t - \langle M^{f,c}, M^c \rangle_t \\ &\quad + \sum_{s \leq t}(f(X_s) - f(X_{s-}))\varphi(X_s, X_{s-})\end{aligned}$$

since both sides are right-continuous in $t \in [0, \zeta[$. Let $\widetilde{K}$ be the purely discontinuous local MAF on $[\![0, \zeta[\![$ with

$$\widetilde{K}_t - \widetilde{K}_{t-} = -f(X_{t-})\varphi(X_{t-}, X_t) - f(X_t)\varphi(X_t, X_{t-}) \qquad \text{for all } t \in [0, \zeta[.$$

Then, by (4.2),

$$\begin{aligned}\Lambda(f*M)_t &= -\tfrac{1}{2}((f*M)_t + (f*M) \circ r_t + f(X_t)\varphi(X_t, X_{t-}) + \widetilde{K}_t) \\ &= \tfrac{1}{2}\left(\int_0^t f(X_{s-})\, dK_s + 2\int_0^t f(X_{s-})\, dA_s + \langle M^{f,c}, M^c \rangle_t \right.\\ &\qquad \left. - \sum_{s \leq t}(f(X_s) - f(X_{s-}))\varphi(X_s, X_{s-}) - \widetilde{K}_t\right).\end{aligned}$$

Thus,

$$\int_0^t f(X_{s-})\, d\Lambda(M)_s$$
$$= \Lambda(f*M)_t - \tfrac{1}{2}\langle M^{f,c}, M^c\rangle_t$$
$$+ \tfrac{1}{2}\int_0^t \int_E (f(y)-f(X_s))\varphi(y,X_s)N(X_s,dy)\,dH_s$$
$$= \tfrac{1}{2}\int_0^t f(X_{s-})\,dK_s + \int_0^t f(X_{s-})\,dA_s$$
$$- \tfrac{1}{2}\sum_{s\le t}(f(X_s)-f(X_{s-}))\varphi(X_s,X_{s-}) - \tfrac{1}{2}\widetilde{K}_t$$
$$+ \tfrac{1}{2}\int_0^t \int_E (f(y)-f(X_s))\varphi(y,X_s)N(X_s,dy)\,dH_s.$$

It now follows from (3.17)–(3.18) that, $\mathbf{P}_m$-a.e.,

$$\int_0^t f(X_{s-})\,d\Lambda(M)_s = \int_0^t f(X_{s-})\,dA_s \qquad \text{for all } t\in[0,\zeta[.$$

This proves the theorem. □

Note that if $f,g\in\mathcal{F}_{\text{loc}}$, then $fg\in\mathcal{F}_{\text{loc}}$.

THEOREM 4.2. *Let $f,g\in\mathcal{F}_{\text{loc}}$ and let $M$ be a locally square-integrable MAF on $[\![0,\zeta[\![$ satisfying (3.1). Then $\mathbf{P}_m$-a.e.,*

$$\text{(4.3)} \quad \int_0^t g(X_{s-})\,d\left(\int_0^s f(X_{r-})\,d\Lambda(M)_r\right)$$
$$= \int_0^t f(X_{s-})g(X_{s-})\,d\Lambda(M)_s \qquad \text{for every } t<\zeta,$$

*whenever all of the integrals involved are well defined.*

PROOF. Let $\varphi:E\times E\to\mathbb{R}$ be a Borel function with $\varphi(x,x)=0$ for $x\in E$ such that, $\mathbf{P}_m$-a.e.,

$$\varphi(X_{t-},X_t) = M_t - M_{t-} \qquad \text{for all } t\in\,]0,\zeta[.$$

Let $K_t$ and $\widetilde{K}_t$ be the purely discontinuous local MAF's on $[\![0,\zeta[\![$ with

$$K_t - K_{t-} = -\varphi(X_{t-},X_t) - \varphi(X_t,X_{t-}) \qquad \text{for } t\in\,]0,\zeta[$$

and

$$\widetilde{K}_t - \widetilde{K}_{t-} = -f(X_{t-})\varphi(X_{t-},X_t) - f(X_t)\varphi(X_t,X_{t-}) \qquad \text{for } t\in\,]0,\zeta[,$$

respectively. The left-hand side of (4.3) is then equal to

$$\int_0^t g(X_{s-})\,d\Lambda(f*M)_s - \tfrac{1}{2}\int_0^t g(X_{s-})\,d\langle M^{f,c}, M^c\rangle_s$$
$$+ \tfrac{1}{2}\int_0^t \int_E g(X_s)(f(y)-f(X_s))\varphi(y,X_s)N(X_s,dy)\,dH_s$$
$$= \Lambda(fg*M)_t - \tfrac{1}{2}\langle M^{g,c}, (f*M)^c\rangle_t$$
$$+ \tfrac{1}{2}\int_0^t \int_E (g(y)-g(X_s))f(y)\varphi(y,X_s)N(X_s,dy)\,dH_s$$
$$- \tfrac{1}{2}\int_0^t g(X_{s-})\,d\langle M^{f,c}, M^c\rangle_s$$
$$+ \tfrac{1}{2}\int_0^t \int_E g(X_s)(f(y)-f(X_s))\varphi(y,X_s)N(X_s,dy)\,dH_s$$
$$= \Lambda(fg*M)_t - \tfrac{1}{2}\langle M^{fg,c}, M^c\rangle_t$$
$$+ \tfrac{1}{2}\int_0^t \int_E (f(y)g(y)-f(X_s)g(X_s))\varphi(y,X_s)N(X_s,dy)\,dH_s$$
$$= \int_0^t f(X_{s-})g(X_{s-})\,d\Lambda(M)_s.$$

This proves the theorem. □

Let $\mathcal{J}$ denote the class of stochastic processes that can be written as the sum of an $(\mathcal{F}_t)$-semimartingale $Y$ and $\Lambda(M)$ for a locally square-integrable MAF $M$ on $[\![0,\zeta[\![$ satisfying the condition of Definition 3.3. The last two theorems imply that the following stochastic integral is well defined for integrators $Z \in \mathcal{J}$.

DEFINITION 4.3. For $f \in \mathcal{F}_{\text{loc}}$ and $Z = Y + \Lambda(M) \in \mathcal{J}$, define on, $[0,\zeta[$,

$$\int_0^t f(X_{s-})\,dZ_s := \int_0^t f(X_{s-})\,dY_s + \int_0^t f(X_{s-})\,d\Lambda(M)_s,$$

whenever the latter stochastic integral is well defined.

To establish Itô's formula, we need the following result.

THEOREM 4.4. *Let $f \in \mathcal{F}_{\text{loc}}$ and let $M$ be a locally square-integrable MAF on $[\![0,\zeta[\![$ such that $\int_0^\cdot f(X_{s-})\,d\Lambda(M)$ is well defined on $[0,\zeta[$. Then, for every $t > 0$, $\mathbf{P}_m$-a.e. on $\{t < \zeta\}$,*

$$(4.4) \quad \int_0^t f(X_{s-})\,d\Lambda(M)_s = \lim_{n\to\infty} \sum_{\ell=0}^{n-1} f(X_{\ell t/n})(\Lambda(M)_{(\ell+1)t/n} - \Lambda(M)_{\ell t/n}).$$

*Here, the convergence is in measure with respect to* $\mathbf{P}_{gm}$ *on* $\{t < \zeta\}$ *for every* $g \in L^1(E; m)$ *with* $0 < g \leq 1$ *m-a.e.*

PROOF. By (3.5), $M_s \circ r_t = M_t \circ r_t - M_{t-s} \circ r_{t-s}$ for all $s < t$. Let $\varphi : E \times E \to \mathbb{R}$ be a Borel function with $\varphi(x, x) = 0$ for $x \in E$ such that $\varphi(X_{t-}, X_t) = M_t - M_{t-}$ for all $t \in [0, \zeta[$. Let $K$ be the purely discontinuous local MAF on $[\![0, \zeta[\![$ with

$$K_t - K_{t-} = -\varphi(X_{t-}, X_t) - \varphi(X_t, X_{t-}) \qquad \text{for } t \in ]0, \zeta[.$$

Then, for each fixed $t > 0$, $\mathbf{P}_m$-a.e. on $\{t < \zeta\}$,

$$\lim_{n \to \infty} \sum_{\ell=0}^{n-1} f(X_{\ell t/n})(\Lambda(M)_{(\ell+1)t/n} - \Lambda(M)_{\ell t/n})$$

$$= -\tfrac{1}{2}(f * M)_t - \tfrac{1}{2}(f * K)_t$$

$$+ \tfrac{1}{2} \lim_{n \to \infty} \sum_{\ell=0}^{n-1} f(X_{\ell t/n})(M_{(\ell+1)t/n} \circ r_{(\ell+1)t/n} - M_{\ell t/n} \circ r_{\ell t/n})$$

$$= -\tfrac{1}{2}(f * M)_t - \tfrac{1}{2}(f * K)_t$$

$$- \tfrac{1}{2} \lim_{n \to \infty} \left[ \sum_{\ell=0}^{n-1} f(X_{(\ell+1)t/n})(M_{(\ell+1)t/n} - M_{\ell t/n}) \right] \circ r_t$$

$$= -\tfrac{1}{2}(f * M)_t - \tfrac{1}{2}(f * K)_t - \tfrac{1}{2}(f * M)_t \circ r_t - \tfrac{1}{2}[M^f, M]_t \circ r_t$$

$$= -\tfrac{1}{2}(f * M)_t - \tfrac{1}{2}(f * K)_t - \tfrac{1}{2}(f * M)_t \circ r_t - \tfrac{1}{2}\langle M^{f,c}, M^c \rangle_t$$

$$- \tfrac{1}{2} \sum_{s \leq t}(f(X_{s-}) - f(X_s))\varphi(X_s, X_{s-})$$

$$= \Lambda(f * M)_t + \tfrac{1}{2}\widetilde{K}_t - \tfrac{1}{2}(f * K)_t - \tfrac{1}{2}\langle M^{f,c}, M^c \rangle_t$$

$$- \tfrac{1}{2} \sum_{s \leq t}(f(X_{s-}) - f(X_s))\varphi(X_s, X_{s-})$$

$$= \int_0^t f(X_{s-}) \, d\Lambda(M)_s,$$

where $\widetilde{K}$ in the penultimate equality is the purely discontinuous local MAF on $[\![0, \zeta[\![$ with $\widetilde{K}_s - \widetilde{K}_{s-} = -f(X_{s-})\varphi(X_{s-}, X_s) - f(X_s)\varphi(X_s, X_{s-})$ for $s \in ]0, \zeta[$. □

REMARK 4.5. (i) Theorem 4.4 immediately implies Theorems 3.10 and 4.1.

(ii) By (3.9),

$$(4.5) \quad \int_0^t f(X_{s-})\,d\Lambda(M)_s = \lim_{n\to\infty} \sum_{\ell=0}^{n-1} f(X_{(\ell+1)t/n})(\Lambda(M)_{(\ell+1)t/n} - \Lambda(M)_{\ell t/n})$$

holds in $P_{gm}$-measure on $\{t < \zeta\}$ for any $g \in L^1(E;m)$ with $0 < g \leq 1$ $m$-a.e. Hence, we could denote this stochastic integral by either $\int_0^t f(X_s)\,d\Lambda(M)_s$ or $\int_0^t f(X_s) \circ d\Lambda(M)_s$. Here, $\int_0^t f(X_s) \circ d\Lambda(M)_s$ is the Fisk–Stratonovich type integral: for $t < \zeta$

$$\begin{aligned}(4.6) \quad &\int_0^t f(X_s) \circ d\Lambda(M)_s \\ &:= \lim_{n\to\infty} \sum_{\ell=0}^{n-1} \frac{f(X_{(\ell+1)t/n}) + f(X_{\ell t/n})}{2}(\Lambda(M)_{(\ell+1)t/n} - \Lambda(M)_{\ell t/n}).\end{aligned}$$

(iii) For any $f \in \mathcal{F}_{\mathrm{loc}}$ and $\mathbf{P}_m$-square-integrable MAF $M$, by way of the Riemann sum approximation (4.4), we can extend the stochastic integral $\int_0^t f(X_{s-})\,d\Lambda(M)_s$ without imposing further conditions. Indeed, let $\{G_\ell\}$ be a nest of finely open Borel sets and $f_\ell \in \mathcal{F}_b$ with $f = f_\ell$ $m$-a.e. on $G_\ell$ [see the explanation for the condition (3.11) in Remark 3.9]. By (4.4), we see $\int_0^t f_n(X_{s-})\,d\Lambda(M)_s = \int_0^t f_\ell(X_{s-})\,d\Lambda(M)_s$ for $t < \tau_{G_n}$ and $n < \ell$. We can then define $\int_0^t f(X_{s-})\,d\Lambda(M)_s = \int_0^t f_n(X_{s-})\,d\Lambda(M)_s$ for $t < \tau_{G_n}$ for each $n \in \mathbb{N}$ and, consequently, for all $t < \zeta$ $\mathbf{P}_m$-a.e. More strongly, for $M \in \overset{\circ}{\mathcal{M}}$ and $f \in \mathcal{F}_{\mathrm{loc}}$, we can define $\int_0^t f(X_{s-})\,d\Lambda(M)_s$ for all $t \in [0,\infty[$ $\mathbf{P}_m$-a.e. Indeed, by Remark 3.9(ii), our stochastic integral $f_n * \Lambda(M)$ for $M \in \overset{\circ}{\mathcal{M}}$ agrees with that defined by Nakao [14] on $[0,\zeta[$ $\mathbf{P}_m$-a.e., while the latter is defined as a CAF of $X$ for all $t \geq 0$. This implies that $\lim_{s\uparrow\zeta}(f_n * \Lambda(M))_s$ exists and is finite $\mathbf{P}_m$-a.e. After we extend our definition of stochastic integral $f_n * \Lambda(M)$ beyond $[0,\zeta[$ by

$$(f_n * \Lambda(M))_t = (f_n * \Lambda(M))_\zeta = \lim_{s\uparrow\zeta}(f_n * \Lambda(M))_s \qquad \text{for } t \geq \zeta,$$

$f_n * \Lambda(M)$ becomes a CAF of $X$ on $[0,\infty[$ $\mathbf{P}_m$-a.e. With this extension for each $n < \ell$, we have $\int_0^t f(X_{s-})\,d\Lambda(M)_s = \int_0^t f_\ell(X_{s-})\,d\Lambda(M)_s$ for $t < \sigma_{E\setminus G_n}$, $\mathbf{P}_m$-a.e. Owing to Lemma 3.1(i) and the existence of the limit $\lim_{t\uparrow\zeta} \int_0^t f(X_{s-})\,d\Lambda(M)_s$ $\mathbf{P}_m$-a.e., we obtain the stochastic integral $\int_0^t f(X_{s-})\,d\Lambda(M)_s$, on $[0,\infty[$, $\mathbf{P}_m$-a.e. for any $f \in \mathcal{F}_{\mathrm{loc}}$ and $M \in \overset{\circ}{\mathcal{M}}$, extending the stochastic integral of Nakao [14].

Remark 4.5(iii) says that the stochastic integral $f * \Lambda(M)_t := \int_0^t f(X_{s-})\,d\Lambda(M)_s$ can be defined for $t \in [0,\infty[$, $\mathbf{P}_m$-a.e., for every $f \in \mathcal{F}_{\mathrm{loc}}$ and $M \in \overset{\circ}{\mathcal{M}}$. We shall refine this statement from $m$-almost every starting point $x \in E$ to quasi-every $x \in E$.

STOCHASTIC CALCULUS FOR SYMMETRIC MARKOV PROCESSES 39LEMMA 4.6. *For $f \in \mathcal{F}_{\mathrm{loc}}$ and $M \in \overset{\circ}{\mathcal{M}}$, the stochastic integral $f * \Lambda(M)_t := \int_0^t f(X_{s-}) \, d\Lambda(M)_s$ can be defined for all $t \in [0, \infty[$, $\mathbf{P}_x$-a.s. for q.e. $x \in E$, in particular, $f * \Lambda(M)$ is a CAF of $X$ on $[0, \infty[$.*

PROOF. Since $f \in \mathcal{F}_{\mathrm{loc}}$, we have $\{f_k \mid k \in \mathbb{N}\} \subset \mathcal{F}_b$ and a nest $\{G_k \mid k \in \mathbb{N}\}$ of finely open Borel sets such that $f = f_k$ q.e. on $G_k$. We know that the stochastic integral $f_k * \Lambda(M)$ is defined $\mathbf{P}_x$-a.s. for q.e. $x \in E$. Let $\Xi_k$ be the defining set admitting an $\mathcal{E}$-polar set for the CAF $f_k * \Lambda(M)$ of zero energy and set

$$\Xi := \left\{ \omega \in \bigcap_{k=1}^{\infty} \Xi_k \;\middle|\; \text{for any } k, \ell \in \mathbb{N} \text{ with } k < \ell, \right.$$

$$\int_0^t f_k(X_{s-}(\omega)) \, d\Lambda(M)_s(\omega)$$

$$\left. = \int_0^t f_\ell(X_{s-}(\omega)) \, d\Lambda(M)_s(\omega) \text{ for } t < \sigma_{E \setminus G_k}(\omega) \right\}.$$

Then, $\mathbf{P}_x(\Xi^c) = 0$, $m$-a.e. $x \in E$. Hence, for each $s > 0$, $\mathbf{P}_x(\theta_s^{-1}(\Xi^c)) = P_s(\mathbf{P}_\cdot(\Xi^c))(x) = 0$ for q.e. $x \in E$. Setting $\widehat{\Xi} := \bigcap_{k=1}^{\infty} \Xi_k \cap \bigcap_{s \in \mathbb{Q}_{++}} \theta_s^{-1}(\Xi)$, we have $\mathbf{P}_x(\widehat{\Xi}) = 1$ for q.e. $x \in E$. For $\omega \in \widehat{\Xi}$ with $t < \sigma_{E \setminus G_k}(\omega)$, we can find small $s_0(= s_0(\omega)) > 0$ such that $t + s_0 < \sigma_{E \setminus G_k}(\omega)$. We then see that $t < \sigma_{E \setminus G_k}(\theta_s \omega)$ for any rational $s \in ]0, s_0[$. Hence, for such $\omega$, we have for $k < \ell$ and any rational $s \in ]0, s_0[$

$$\int_s^{t+s} f_k(X_{v-}(\omega)) \, d\Lambda(M)_v(\omega) = \int_s^{t+s} f_\ell(X_{v-}(\omega)) \, d\Lambda(M)_v(\omega).$$

Letting $s \to 0$ and noting that $\omega \in \Xi_k$, $k \in \mathbb{N}$, we have that for $k < \ell$, $f_k * \Lambda(M)_t = f_\ell * \Lambda(M)_t$ for $t < \sigma_{E \setminus G_k}$, $\mathbf{P}_x$-a.s. for q.e. $x \in E$. By Lemma 3.1(i), we know that $\mathbf{P}_x(\lim_{k \to \infty} \sigma_{E \setminus G_k} = \infty) = 1$ for q.e. $x \in E$. Therefore, we obtain that the stochastic integral $f * \Lambda(M)$ defined as in Remark 4.5(4.5) can be established $\mathbf{P}_x$-a.s. for q.e. $x \in E$. This completes the proof. □

THEOREM 4.7 (Generalized Itô formula). *Suppose that $\Phi \in C^2(\mathbb{R}^d)$ and $u = (u_1, \ldots, u_d) \in \mathcal{F}^d$. Then, for q.e. $x \in E$, $\mathbf{P}_x$-a.s. for all $t \in [0, \infty[$,*

$$\Phi(u(X_t)) - \Phi(u(X_0))$$

$$= \sum_{k=1}^d \int_0^t \frac{\partial \Phi}{\partial x_k}(u(X_{s-})) \, du_k(X_s)$$

(4.7)
$$+ \frac{1}{2} \sum_{i,j=1}^d \int_0^t \frac{\partial^2 \Phi}{\partial x_i \, \partial x_j}(u(X_{s-})) \, d\langle M^{u_i, c}, M^{u_j, c} \rangle_s$$

$$+ \sum_{s \leq t} \Bigg( \Phi(u(X_s)) - \Phi(u(X_{s-}))$$
$$- \sum_{k=1}^{d} \frac{\partial \Phi}{\partial x_k}(u(X_{s-}))(u_k(X_s) - u_k(X_{s-})) \Bigg).$$

PROOF. Note that both sides appearing in (4.7) are $\mathbf{P}_x$-a.s. defined for q.e. $x \in E$ in view of Lemma 4.6. First, we show this Itô formula (4.7) under $\mathbf{P}_m$ for a fixed $t \geq 0$. Note that $\Phi \circ u \in \mathcal{F}_{\mathrm{loc}}$ and that

$$u_k(X_t) = u_k(X_0) + M_t^{u_k} + N_t^{u_k} = u_k(X_0) + M_t^{u_k} + \Lambda(M^{u_k})_t.$$

This version of Itô's formula follows from Theorems 3.7 and 4.4 by a line of reasoning similar to that used to prove Itô's formula for semimartingales (cf. [9]). Since both sides in (4.7) are right-continuous, (4.7) holds under $\mathbf{P}_m$.

Second, we refine the starting point. Recall that $\Omega$ consists of rcll paths. Let $I_t(\omega)$ be the difference of the left-hand side and the right-hand side of (4.7). Let $\Xi$ be the intersection of all of the defining sets of AF's appearing in the formula and $\{\omega \in \Omega \mid I_t(\omega) = 0, \forall t \in [0, \infty[\}$. Then, $\mathbf{P}_x(\Xi^c) = 0$, $m$-a.e. $x \in E$. Let $\widehat{\Xi}$ be the intersection of the defining sets of AF's appearing in the formula and $\bigcap_{s \in \mathbb{Q}_{++}} \theta_s^{-1}(\Xi)$. We then have $\mathbf{P}_x(\widehat{\Xi}) = 1$ for q.e. $x \in E$, as in the proof of Lemma 4.6. Take $\omega \in \widehat{\Xi}$. For any positive rational $s > 0$, we then have $I_t(\theta_s \omega) = 0$, that is,

$$\Phi(u(X_{t+s}(\omega))) - \Phi(u(X_s(\omega)))$$
$$= \sum_{k=1}^{d} \int_s^{t+s} \frac{\partial \Phi}{\partial x_k}(u(X_{v-}(\omega))) \, du_k(X_v(\omega))$$
$$+ \frac{1}{2} \sum_{i,j=1}^{d} \int_s^{t+s} \frac{\partial^2 \Phi}{\partial x_i \partial x_j}(u(X_{v-}(\omega))) \, d\langle M^{u_i,c}, M^{u_j,c} \rangle_v(\omega)$$
$$+ \sum_{s < v \leq t+s} \Bigg( \Phi(u(X_v(\omega))) - \Phi(u(X_{v-}(\omega)))$$
$$- \sum_{k=1}^{d} \frac{\partial \Phi}{\partial x_k}(u(X_{v-}(\omega)))(u_k(X_v(\omega)) - u_k(X_{v-}(\omega))) \Bigg).$$

Letting $s \to 0$ and using the right-continuity of $s \mapsto u(X_s)$ and stochastic integrals, we have $I_t(\omega) = 0$. This completes the proof. $\square$

Z.-Q. Chen
Department of Mathematics
University of Washington
Seattle, Washington 98195
USA
E-mail: zchen@math.washington.edu

K. Kuwae
Department of Mathematics
Faculty of Education
Kumamoto University
Kumamoto 860-8555
Japan
E-mail: kuwae@gpo.kumamoto-u.ac.jp

P. J. Fitzsimmons
Department of Mathematics
University of California at San Diego
La Jolla, California 92093-0112
USA
E-mail: pfitzsim@ucsd.edu

T.-S. Zhang
School of Mathematics
University of Manchester
Sackville Street, Manchester M60 1QD
United Kingdom
E-mail: tzhang@maths.manchester.ac.uk